\documentclass[reqno]{amsart}

\usepackage{general,resistance,ifthen,txfonts,text_shortcuts}
\usepackage{cancel,xy}

\usepackage[bookmarks, colorlinks=true, pdfstartview=FitV, linkcolor=blue, citecolor=blue, urlcolor=blue]{hyperref}
\usepackage{cite}  

\newcommand{\clean}{true}  
\newcommand{\version}[2]{\ifthenelse{\equal{\clean}{true}}{#1}{{\footnotesize #2}}}
\newcommand{\linenopax}{} 

\numberwithin{equation}{section} \numberwithin{theorem}{section}

\begin{document}

\title{A {H}ilbert space approach to effective resistance metric}

\author{Palle E. T. Jorgensen}
\address{University of Iowa, Iowa City, IA 52246-1419 USA}
\email{jorgen@math.uiowa.edu}

\author{Erin P. J. Pearse}
\address{University of Iowa, Iowa City, IA 52246-1419 USA}
\email{epearse@math.uiowa.edu}

\thanks{The work of PETJ was partially supported by NSF grant DMS-0457581. The work of EPJP was partially supported by the University of Iowa Department of Mathematics NSF VIGRE grant DMS-0602242.}

\begin{abstract}
A resistance network is a connected graph $(G,c)$. The conductance function $c_{xy}$ weights the edges, which are then interpreted as conductors of possibly varying strengths. The Dirichlet energy form $\mathcal E$ produces a Hilbert space structure (which we call the energy space ${\mathcal H}_{\mathcal E}$) on the space of functions of finite energy.  

We use the reproducing kernel $\{v_x\}$ constructed in \cite{DGG} to analyze the effective resistance $R$, which is a natural metric for such a network. It is known that when $(G,c)$ supports nonconstant harmonic functions of finite energy, the effective resistance metric is not unique. The two most natural choices for $R(x,y)$ are the ``free resistance'' $R^F$, and the ``wired resistance'' $R^W$. We define $R^F$ and $R^W$ in terms of the functions $v_x$ (and certain projections of them). This provides a way to express $R^F$ and $R^W$ as norms of certain operators, and explain $R^F \neq R^W$ in terms of Neumann vs. Dirichlet boundary conditions. We show that the metric space $(G,R^F)$ embeds isometrically into ${\mathcal H}_{\mathcal E}$, and the metric space $(G,R^W)$ embeds isometrically into the closure of the space of finitely supported functions; a subspace of ${\mathcal H}_{\mathcal E}$.

Typically, $R^F$ and $R^W$ are computed as limits of restrictions to finite subnetworks. A third formulation $R^{tr}$ is given in terms of the trace of the Dirichlet form $\mathcal E$ to finite subnetworks. A probabilistic approach shows that in the limit, $R^{tr}$ coincides with $R^F$. This suggests a comparison between the probabilistic interpretations of $R^F$ vs. $R^W$.

\end{abstract}




  \date{\bf\today. \q \version{}{Rough version (for editing).}}

\maketitle

\setcounter{tocdepth}{1}
{\small \tableofcontents}

\allowdisplaybreaks


This paper concerns an analysis of the resistance metric (also called ``effective resistance'') on infinite networks, with emphasis on the role of functions of finite energy. These are real- or complex-valued functions on the set of vertices of the network, and their energy is computed via a Dirichlet form (denoted \energy) which takes into account the weights (conductances) of the edges of the network. As discussed in \cite{DGG}, the Dirichlet energy form \energy gives a natural Hilbert space structure to the set of finite-energy functions; we call this the ``energy space'' and denote it \HE. The close relationship between \energy and the (typically unbounded) network Laplacian \Lap is also developed in \cite{DGG}. In the absence of $L^2$ conditions, this relationship is more subtle than described in the classical theory of quadratic forms and self-adjoint operators, as given by Kato, for example.  

The Hilbert space \HE does not come equipped with a natural o.n.b., but it does carry a natural reproducing kernel, which is indexed by the vertices $x$ of the network, and is denoted $\{v_x\}$. We make extensive use of this ``energy kernel'' and its properties as developed in \cite{DGG}. Precise definitions are given in \S\ref{sec:The-energy-space}. 

Effective resistance has been previously studied in a Hilbert space context; however, our approach is novel in that it is completely intrinsic to the Hilbert space \HE. For example, there is no reference to space of functions defined on the edges of \Graph as in \cite{Lyons:ProbOnTrees}, and there is no use of potential-theoretic methods, as in \cite{Kig03}. Our theory is developed via a reproducing kernel which springs forth directly from Riesz's lemma, and allows us to express effective resistance in terms of operators. In particular, our approach allows one to treat the effective resistance metrics globally, i.e., directly in terms of Hilbert space norms and independent of any limits taken with respect to finite subnetworks. 
 
\pgap

\S\ref{sec:introduction} contains a brief encapsulation of the results of \cite{DGG} which will be necessary for the current study. In particular, the orthogonal decomposition of \HE into finitely-supported and harmonic functions, and the corresponding projections. 
Dipoles are finite-energy functions which are harmonic except at two vertices. The energy kernel consists of dipoles, and are important for computing the resistance metric.

\S\ref{sec:Resistance-metric-on-finite-networks} introduces the effective resistance as a metric on finite networks and gives several equivalent formulations. We also discuss how certain infinite networks may give rise to nonuniqueness of solutions with regard to these formulations, and why this therefore leads to different notions of resistance metric in infinite networks.
\S \ref{sec:Resistance-metric-on-infinite-networks} discusses the two most well-known resistance metrics on infinite networks: ``free resistance'' $R^F$ and ``wired resistance'' $R^W$ (in the terminology of \cite{Lyons:ProbOnTrees}; the respective terms ``limit'' and ``minimal'' are also common in the older literature). We give formulas for the free and wired resistances in parallel to the list of equivalent formulations developed in the previous section for finite networks. The Hilbert space structure of \HE (i.e., certain projections and operator norms) gives a clear explanation of why $R^F(x,y) \geq R^W(x,y)$ in general, and how harmonic functions can produce a strict inequality. We also consider $R^F$ vs. $R^W$ in terms of Neumann vs. Dirichlet boundary conditions and probabilistic interpretations. The harmonic resistance $R^\hrm := R^F - R^W$ is introduced in \S\ref{sec:Harmonic-resistance} and related to the boundary resistance $R^\bdy$ (roughly, the ``voltage drop at \iy''). Neither $R^\hrm$ nor $R^\bdy$ are metrics, in contrast to $R^F$ and $R^W$. In Remark~\ref{rem:comparison-to-resistance-forms}, we give a comparison with the theory of resistance forms of \cite{Kig01,Kig03,Kig09}.

\S\ref{sec:trace-resistance} introduces a third resistance metric on infinite networks; like the others, it is computed as a limit of restrictions to finite subnetworks. The ``trace resistance'' $R^{tr}$ is given in terms of the trace of \energy to finite subnetworks. A probabilistic approach shows that in the limit, $R^{tr}$ coincides with $R^F$. 

\S\ref{sec:Comparison-to-other-metrics} compares the resistance metric(s) with related notions of distance on a network, including the geodesic distance. We also consider the effective resistance between two probability measures as a natural generalization of our earlier formulations for $R^F$ and $R^W$.

The work of von Neumann and Schoenberg gives conditions under which a metric space embeds in a Hilbert space (so that the original metric is recovered in terms of a normed difference of the embedded points). \S\ref{sec:Construction-of-HE} shows that these conditions are satisfied by the effective resistances discussed above. Moreover, up to unitary isomorphism, the embedding sends the metric space $(G,R^F)$ into \HE, and the metric space $(G,R^W)$ into \Fin. Under this embedding, the vertex $x$ is sent to the element $v_x$ of the energy kernel.

In \S\ref{sec:examples}, we give some elementary but illuminating examples.

\pgap

Discrete potential theory and its relation to random walks on graphs is an old and well-studied subject (for trees and Cayley graphs of groups in particular) and we will not attempt to give complete references. Three excellent and fairly comprehensive treatments are \cite{Soardi94}, \cite{Woess00} and \cite{Lyons:ProbOnTrees}. We also recommend \cite{DoSn84, LevPerWil08, Peres99} for introductory material and \cite{TerryLyons, Car73a, Woess00}, and the foundational paper \cite{Nash-Will59} for more specific background. With regard to infinite graphs and finite-energy functions, see \cite{Soardi94, SoardiWoess91, CaW92, Dod06, PicWoess90, PicWoess88, Wo86, Thomassen90} and also some recent work by Georgakopoulos (see the arXiv).

Effective resistance and resistance metric are studied extensively in \cite{Kig03,Kig09}, and also in \cite{Kig01}. Connections between shorting and Schur complement (trace) are studied in \cite{Metz}.

\subsection*{Acknowledgements}
The authors are grateful to Jun Kigami, Peter M\"{o}rters, Elmar Teufl, and Wolfgang Woess for helpful conversations, suggestions, and answers to our questions.

\section{Introduction}
\label{sec:introduction}

\subsection{Basic terms}
\label{sec:electrical-resistance-networks}

We now proceed to introduce the key notions used throughout this paper: resistance networks, the energy form \energy, the Laplace operator \Lap, and the elementary relations amongst them.

\begin{defn}\label{def:ERN}
  A \emph{resistance network} is a connected graph $(\Graph,\cond)$, where \Graph is a graph with vertex set \verts, and $\cond:\verts \times \verts \to \bR^+$ is the \emph{conductance function} which defines adjacency by $x \nbr y$ iff $c_{xy}>0$, for $x,y \in \verts$. Conductance is symmetric and nonnegative:  $\cond_{xy} = \cond_{yx} \in [0,\iy)$. The total conductance at a vertex is written $\cond(x) := \sum_{y \nbr x} \cond_{xy}$, and it is required that $\cond(x) < \iy$. Note that we allow vertices of infinite degree and that $\cond(x)$ need not be a bounded function on \verts. The notation \cond may be used to indicate the multiplication operator $(\cond v)(x) := \cond(x) v(x)$.  
  \version{}{\marginpar{Is it ok to say ``basis'' here?}}
\end{defn}

In Definition~\ref{def:ERN}, ``connected'' means simply that for any $x,y \in \verts$, there is a finite sequence $\{x_i\}_{i=0}^n$ with $x=x_0$, $y=x_n$, and $\cond_{x_{i-1} x_i} > 0$, $i=1,\dots,n$. Conductance is the reciprocal of resistance, so one can think of $(\Graph,\cond)$ as a network of nodes \verts connected by resistors of resistance $\cond_{xy}^{-1}$. We may assume there is at most one edge from $x$ to $y$, as two conductors $\cond^1_{xy}$ and $\cond^2_{xy}$ connected in parallel can be replaced by a single conductor with conductance $\cond_{xy} = \cond^1_{xy} + \cond^2_{xy}$. Also, we assume $\cond_{xx}=0$ so that no vertex has a loop, as electric current will never flow along a conductor connecting a node to itself.\footnote{Nonetheless, self-loops may be useful for technical considerations: one can remove the periodicity of a random walk by allowing self-loops. This can allow one to obtain a ``lazy walk'' which is ergodic, and hence more tractable. See, for example, \cite{LevPerWil08, Lyons:ProbOnTrees}.}

\begin{defn}\label{def:graph-laplacian}
  The \emph{Laplacian} on \Graph is the linear difference operator 
  which acts on a function $v:\verts \to \bC$ by
  \linenopax
  \begin{equation}\label{eqn:def:laplacian}
    (\Lap v)(x) :
    = \sum_{y \nbr x} \cond_{xy}(v(x)-v(y)).
  \end{equation}
  A \fn $v:\verts \to \bC$ is \emph{harmonic} iff $\Lap v(x)=0$ for each $x \in \verts$. 
\end{defn}

We have adopted the physics convention (so that the spectrum is nonnegative) and thus our Laplacian is the negative of the one commonly found in the PDE literature; e.g., \cite{Kig01}, \cite{Str06}. The network Laplacian \eqref{eqn:def:laplacian} should not be confused with the renormalized Laplace operator $\cond^{-1/2} \Lap \cond^{-1/2}$ which appears in the literature on spectral graph theory (e.g., \cite{Chu01}).

\begin{defn}\label{def:Prob-operator}
  The \emph{(probabilistic) transition operator} is defined pointwise for functions on \verts by
  \linenopax
  \begin{align}\label{eqn:def:Prob-trans-oper}
    \Prob u(x) = \sum_{y \nbr x} p(x,y) u(y), 
    \q\text{for } 
    p(x,y) = \frac{\cond_{xy}}{\cond(x)},
  \end{align}
  so that $\Lap = \cond(\id - \Prob)$. Note that the harmonic functions are precisely the fixed points of \Prob, and $v=u+k\one$ for $k \in \bC$ implies that $\Prob v = \Prob u + k \one$, so \eqref{eqn:def:Prob-trans-oper} is independent of the representative chosen for $u$.

  The function $p(x,y)$ gives transition probabilities, i.e., the probability that a random walker currently at $x$ will move to $y$ with the next step. Since
  \linenopax
  \begin{align}\label{eqn:def:reversible}
    \cond(x) p(x,y) = \cond(y) p(y,x),
  \end{align}  
  the transition operator \Prob determines a \emph{reversible} Markov process with state space \verts; see \cite{DoSn84,LevPerWil08,Lyons:ProbOnTrees,Peres99}. 
\end{defn}

\begin{defn}\label{def:exhaustion-of-G}
  An \emph{exhaustion} of \Graph is an increasing sequence of finite and connected subgraphs $\{\Graph_k\}$, so that $\Graph_k \ci \Graph_{k+1}$ and $\Graph = \bigcup \Graph_k$.
  \glossary{name={$\{\Graph_k\}$},description={exhaustion of a network},sort=G,format=textbf}
\end{defn}

\begin{defn}\label{def:infinite-vertex-sum}  
  The notation
  \linenopax
  \begin{equation}\label{eqn:def:infinite-sum}
    \sum_{x \in \verts} := \lim_{k \to \iy} \sum_{x \in \Graph_k}
  \end{equation}
  is used whenever the limit is independent of the choice of exhaustion $\{\Graph_k\}$ of \Graph. This is clearly justified, for example, whenever the sum has only finitely many nonzero terms, or is absolutely convergent as in the definition of \energy just below.
\end{defn}

\begin{defn}\label{def:graph-energy}
  The \emph{energy} of functions $u,v:\verts \to \bC$ is given by the (closed, bilinear) Dirichlet form
  \linenopax
  \begin{align}\label{eqn:def:energy-form}
    \energy(u,v)
    := \frac12 \sum_{x \in \verts}  \sum_{y \in \verts} \cond_{xy} \left(\cj{u(x)}-\cj{u(y)}\right)(v(x)-v(y)),
  \end{align}
  with the energy of $u$ given by $\energy(u) := \energy(u,u)$.
  The \emph{domain of the energy} is
  \linenopax
  \begin{equation}\label{eqn:def:energy-domain}
    \dom \energy = \{u:\verts \to \bC \suth \energy(u)<\iy\}.
  \end{equation}
\end{defn}

Since $\cond_{xy}=\cond_{yx}$ and $\cond_{xy}=0$ for nonadjacent vertices, the initial factor of $\frac12$ in \eqref{eqn:def:energy-form} implies there is exactly one term in the sum for each edge in the network. 

\subsection{The energy space \HE}
\label{sec:The-energy-space}

For remainder of this paper, let $o$ denote a fixed reference vertex (the ``origin''). It will be immediate that all results are independent of the choice of $o$. Note from Definition~\ref{def:graph-energy} that $\energy(u)=0$ iff $u$ is constant. Let \one denote the constant function with value 1, so that $\ker \energy = \bC \one$. 

\begin{defn}\label{def:H_energy}\label{def:The-energy-Hilbert-space}
  The energy form \energy is symmetric and positive definite on $\dom \energy$. Then $\dom \energy / \bC \one$ is a vector space with inner product and corresponding norm given by
  \linenopax
  \begin{equation}\label{eqn:energy-inner-product}
    \la u, v \ra_\energy := \energy(u,v)
    \q\text{and}\q
    \|u\|_\energy := \energy(u,u)^{1/2}.
  \end{equation}
  The \emph{energy Hilbert space} \HE is $\dom \energy / \bC \one$. 
\end{defn}

\begin{defn}\label{def:vx}\label{def:energy-kernel}
  Let $v_x$ be defined to be the unique element of \HE for which
  \linenopax
  \begin{equation}\label{eqn:def:vx}
    \la v_x, u\ra_\energy = u(x)-u(o),
    \qq \text{for every } u \in \HE.
  \end{equation}
  The collection $\{v_x\}_{x \in \verts}$ forms a reproducing kernel for \HE (\cite[Cor.~2.7]{DGG}); we call it the \emph{energy kernel} and \eqref{eqn:def:vx} implies that its span is dense in \HE. In \cite[Prop.~4.3]{Kig03}, it is stated that this family is a reproducing kernel for the form \energy.
  Functions with the property \eqref{eqn:def:vx} have also appeared in \cite{Metz,BakerRumely,MuYaYo}, and \cite[Def.~4.1]{Kig03}. 
\end{defn}

  Note that $v_o$ corresponds to a constant function, since $\la v_o, u\ra_\energy = 0$ for every $u \in \HE$. Since $v_o \in \ker \energy$, it may be ignored or omitted from any set or sum indexed by $\{x \in \verts\}$.
    
\begin{defn}\label{def:dipole}
  A \emph{dipole} is any $v \in \HE$ satisfying the pointwise identity $\Lap v = \gd_x - \gd_y$ for some vertices $x,y \in \verts$. One can check that $\Lap v_x = \gd_x - \gd_o$; cf. \cite[Lemma~2.13]{DGG}.
\end{defn}

\begin{defn}\label{def:Fin}
  For $v \in \HE$, one says that $v$ has \emph{finite support} iff there is a finite set $F \ci \verts$ such that $v(x) = k \in \bC$ for all $x \notin F$. That is, the set of functions of finite support in \HE is 
  \linenopax
  \begin{equation}\label{eqn:span(dx)}
    \spn\{\gd_x\} = \{u \in \dom \energy \suth u(x)=k \text{ for some $k$, for all but finitely many } x \in \verts\},
  \end{equation}
  where $k$ is some constant depending only on $u$, and $\gd_x$ is the Dirac mass at $x$, i.e., the element of \HE containing the characteristic function of the singleton $\{x\}$. It is immediate from \eqref{eqn:def:energy-form} that $\energy(\gd_x) = \cond(x)$, whence $\gd_x \in \HE$.
  Define \Fin to be the closure of $\spn\{\gd_x\}$ with respect to \energy. 
\end{defn}

\begin{defn}\label{def:Harm}
  The set of harmonic functions of finite energy is denoted
  \linenopax
  \begin{equation}\label{eqn:Harm}
    \Harm := \{v \in \HE \suth \Lap v(x) = 0, \text{ for all } x \in \verts\}.
  \end{equation}
  Note that this is independent of choice of representative for $v$ in virtue of \eqref{eqn:def:laplacian}.
\end{defn}

\begin{lemma}[{\cite[2.11]{DGG}}]
  \label{thm:<delta_x,v>=Lapv(x)}
  For any $x \in \verts$, one has $\la \gd_x, u \ra_\energy = \Lap u(x)$.
\end{lemma}

The following result follows easily from Lemma~\ref{thm:<delta_x,v>=Lapv(x)}; cf.~\cite[Thm.~2.15]{DGG}.

\begin{theorem}[Royden decomposition]
  \label{thm:HE=Fin+Harm}
  $\HE = \Fin \oplus \Harm$.
\end{theorem}

\begin{remark}[Reproducing kernels for \Fin and \Harm]
  \label{rem:3-repkernels}
  Throughout the sequel, we use the notation \Pfin for orthogonal projection to \Fin and \Phar for orthogonal projection to \Harm. Also, for an element $v_x$ of the energy kernel, we write $f_x := \Pfin v_x$ and $h_x := \Phar v_x$. The characteristic property of reproducing kernels behaves well with respect to orthogonal projection, and so $\{f_x\}_{x \in \Graph}$ is a reproducing kernel for \Fin and $\{h_x\}_{x \in \Graph}$ is a reproducing kernel for \Harm. In particular, $\spn\{h_x\}_{x \in \Graph}$ is a dense subspace of \Harm.
\end{remark}

\begin{remark}[Real-valued vs. complex-valued functions]
  \label{rem:real-vs-complex}
  The setting laid out in this section is valid for \bC-valued functions, and indeed, the associated spectral theory of \Lap requires this. However, \bC-valued functions will not be necessary for the purposes of this paper. Thus, from this point on, the discussion will concern only \bR-valued functions, as this will simplify the exposition.  Even when this simplifying assumption is not made, \cite[Lem.~2.24]{DGG} shows that for the kernel elements $v_x$, $f_x$, and $h_x$ discussed in the previous remark, one can always choose a representative  which is \bR-valued.
\end{remark}

\section{Effective resistance}
\label{sec:Effective-resistance}

Our main concern is the metric properties of effective resistance on infinite networks. Infinite networks may offer several distinct such metrics, each one reflecting separate dynamical or potential-theoretic features and conclusions for $(\Graph,\cond)$. Adopting the terminology and notation of \cite{Lyons:ProbOnTrees}, we consider the free resistance $R^F$ and the wired resistance $R^W$. Later, in \S\ref{sec:trace-resistance}, we introduce and study the trace resistance $R^\trc$. For an infinite network, each of these is initially defined in terms of a limit of finite subnetworks associated to an exhaustion of $(\Graph,\cond)$. However, the introduction of Hilbert space allows one to treat these metrics globally, i.e., directly with reference to a norm in Hilbert space and independent of any limits taken over finite subnetworks. This is carried out for the free resistance $R^F$ in Theorem~\ref{thm:free-resistance} and for the wired resistance $R^W$ in Theorem~\ref{thm:wired-resistance}. 

For the sections to follow, it will be helpful to have the following terminology about currents, which play the role of a discrete vector field. The \emph{drop operator} \drp is an isometry relating functions on the vertices of $(\Graph,\cond)$ to functions on edges.

\begin{defn}\label{def:currents}
  A \emph{current} is an antisymmetric function $\curr: \verts \times \verts \to \bR$ which is supported on the support of \cond. For such functions \curr and $J$, we have
  \linenopax
  \begin{equation}\label{eqn:def:dissipation}
    \diss(\curr,J) := \frac12 \sum_{(x,y) \in \edges} \cond_{xy}^{-1} \curr(x,y) J(x,y),
  \end{equation}
  the \emph{dissipation} of a current is $\diss(\curr) := \diss(\curr,\curr)$. The inner product $\la \curr,J\ra_\diss := \diss(\curr,J)$ makes $\dom \diss := \{\curr \suth \diss(\curr)<\iy\}$ into a Hilbert space which we call the \emph{dissipation space} \HD.
  
  For $a,z \in \verts$, we say \curr is a \emph{current flow from $a$ to $z$} and write $\curr \in \Flo(a,z)$ iff it satisfies
  \linenopax
  \begin{equation}\label{eqn:Kirchhoff-nonhomog}
    \sum_{y \nbr x} \curr(x,y)
    = \gd_a - \gd_z,
  \end{equation}
  a nonhomogeneous version of Kirchhoff's law. 
  If $u \in \HE$, then the \emph{induced current} is defined by Ohm's law:
  \linenopax
  \begin{align}\label{eqn:Ohms-law}
    \curr_u(x,y) = \cond_{xy}(u(x) - u(y)).
  \end{align}
\end{defn}
  Note that $\drp:\HE \to \HD$ by $\drp u = \curr_u$ is an isometry, and that \curr minimizes \diss over $\Flo(x,y)$ if and only if $\curr = \drp u$ for $u$ which minimizes \energy over $\{v \in \dom \energy \suth \Lap v = \gd_x - \gd_y\}$; see \cite[\S3 and \S10]{OTERN} for details.

\subsection{Resistance metric on finite networks}
\label{sec:Resistance-metric-on-finite-networks}

We make the standing assumption that the network is finite in \S\ref{sec:Resistance-metric-on-finite-networks}. However, the results actually remain true on any network for which $\Harm = 0$.

\begin{defn}\label{def:effective-resistance}
  If a current of one amp is inserted into the \ERN at $x$ and withdrawn at $y$, then the \emph{(effective) resistance} $R(x,y)$ is the voltage drop between the vertices $x$ and $y$.
  \glossary{name={$R$},description={resistance metric},sort=r,format=textbf}
\end{defn}

\begin{theorem}\label{thm:effective-resistance-metric}
  The resistance $R(x,y)$ has the following equivalent formulations:
  \linenopax
  \begin{align}
    R(x,y)
    &= \{v(x)-v(y) \suth \Lap v = \gd_x-\gd_y\}
    \label{eqn:def:R(x,y)-Lap} \\
    &= \{\energy(v) \suth \Lap v = \gd_x-\gd_y\} \label{eqn:def:R(x,y)-energy} \\
    &= \min \{\diss(\curr) \suth \curr \in \Flo(x,y)\} \label{eqn:def:R(x,y)-diss} \\
    &= 1/\min \{\energy(v) \suth v(x)=1, v(y)=0, v \in \dom\energy\} \label{eqn:def:R(x,y)-R} \\
    &= \min \{\gk \geq 0 \suth |v(x)-v(y)|^2 \leq \gk \energy(v), v \in \dom\energy\} \label{eqn:def:R(x,y)-S} \\
    &= \sup \{|v(x)-v(y)|^2 \suth \energy(v) \leq 1, v \in \dom\energy\} \label{eqn:def:R(x,y)-sup}.
  \end{align}
\end{theorem}

We leave the proof of Theorem~\ref{thm:effective-resistance-metric} as an exercise; we suggest using the energy kernel to take several shortcuts. A complete proof appears in \cite[Thm.~5.2]{OTERN} and patches some holes in the literature. The authors first learned of the effective resistance metric from \cite{Pow76b}) and \cite{Kig01,Kig03,Kig09,Str06}, respectively; we have not seen \eqref{eqn:def:R(x,y)-energy} in the literature previously. Taking the minimum (rather than the infimum) in \eqref{eqn:def:R(x,y)-diss}, etc, is justified because a quadratic form always attains its minimum on a closed convex set. Effective resistance is defined in \cite[\S8]{Peres99} as the ratio $(v(x)-v(y))/\sum_{z \nbr x} \cond_{xz} (v(x)-v(z))$; our formulation corresponds to normalizing the current flow so that the denominator is 1.

\begin{figure}
    \centering
    \includegraphics{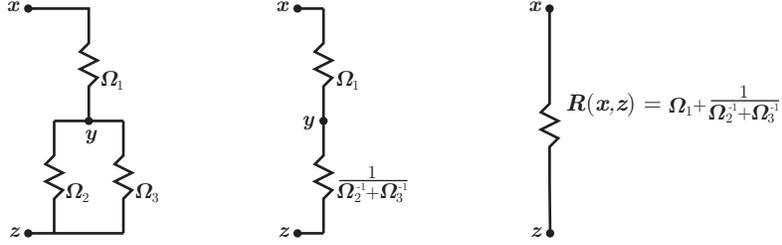}
    \caption{\captionsize Effective resistance as network reduction to a trivial network. This basic example uses parallel reduction followed by series reduction; see Remark~\ref{rem:Resistance-distance-via-network-reduction}.}
    \label{fig:network-reduction}
  \end{figure}

\begin{remark}[Resistance distance via network reduction]
  \label{rem:Resistance-distance-via-network-reduction}
  Let $H$ be a (connected) planar subnetwork of a finite network $G$ and pick any $x,y \in H$. Then $H$ may be reduced to a trivial network consisting only of these two vertices and a single edge between them via the use of three basic transformations: (i) series reduction, (ii) parallel reduction, and (iii) the $\nabla$-\textsf{Y} transform \cite{Epifanov66,Truemper89}. Each of these transformations preserves the resistance properties of the subnetwork, that is, for $x,y \in G \less H$, $R(x,y)$ remains unchanged when these transformations are applied to $H$. The effective resistance between $x$ and $y$ may be interpreted as the resistance of the resulting single edge. An elementary example is shown in Figure~\ref{fig:network-reduction}. A more sophisticated technique of network reduction is given by the trace (Schur complement) construction of Remark~\ref{rem:Resistance-distance-via-Schur-complement}, which subsumes (i) and (iii).
\end{remark}

We record the following simple fact for future reference. 

\begin{lemma}\label{thm:max-and-min-of-vx}
  If $v \in \HE$ is a dipole on a finite network with $\Lap v = \gd_x - \gd_y$, then $v$ has its maximum at $x$ and minimum at $y$.
  \begin{proof}
    This follows by the minimum principle for harmonic functions on the finite subnetwork $\verts\less\{x,y\}$; cf.~\cite[\S2.1]{Lyons:ProbOnTrees} or \cite{LevPerWil08}, for example. 
  \end{proof}
\end{lemma}

The following result is well-known (see, e.g. \cite[\S2.3]{Kig01}), but the proof given here is substantially simpler than most others found in the literature.

\begin{lemma}\label{thm:R-is-a-metric}
  $R$ is a metric.
  \begin{proof}
    Symmetry and positive definiteness are immediate from \eqref{eqn:def:R(x,y)-energy}, we use \eqref{eqn:def:R(x,y)-Lap} and the energy kernel to check the triangle inequality. Let $v_1 = v_x - v_y$, $v_2 = v_y - v_z$, and $v_3 := v_1 + v_2$. Then
    \linenopax
    \begin{align*}
      R(x,z)
      = v_3(x)-v_3(z)
      &= v_1(x)-v_1(z) + v_2(x)-v_2(z) \\
      &\leq v_1(x)-v_1(y) + v_2(y)-v_2(z)
      = R(x,y) + R(y,z),
    \end{align*}
    because $y$ is the minimum of $v_1$ and the maximum of $v_2$ by Lemma~\ref {thm:max-and-min-of-vx}.
  \end{proof}
\end{lemma}

\subsection{Resistance metric on infinite networks}
\label{sec:Resistance-metric-on-infinite-networks}

There are challenges in the extension of the results of the previous section to infinite networks. The existence of nonconstant harmonic functions $h \in \dom \energy$ implies the nonuniqueness of solutions to $\Lap u = f$ in \HE, and hence \eqref{eqn:def:R(x,y)-Lap}--\eqref{eqn:def:R(x,y)-diss} are no longer well-defined. Even worse, the two most common formulations in the literature, \eqref{eqn:def:R(x,y)-R} and \eqref{eqn:def:R(x,y)-S}, do not remain equivalent for networks with $\Harm \neq 0$. Explaining how and why this disparity can occur on infinite networks comprises a large part of the motivation for this paper.

Two natural choices for extension lead to the free resistance $R^F$ and the wired resistance $R^W$. In this section, we examine the cause of a strict inequality $R^W(x,y) < R^F(x,y)$. 
\begin{enumerate}
  \item Theorem~\ref{thm:free-resistance} shows how $R^F$ corresponds to choosing solutions to $\Lap u = \gd_x - \gd_y$ from the energy kernel, and how it corresponds to currents which are decomposable in terms of paths. In \S\ref{sec:probabilistic-interp-of-v_x-and-f_x}, the latter leads to a probabilistic interpretation which provides for a relation to the trace of the resistance discussed in \S\ref{sec:trace-resistance}. See also Corollary~\ref{thm:trace-resistance-is-free-resistance}.
  \item Theorem~\ref{thm:wired-resistance} shows how $R^W$ corresponds to solutions obtained by projecting elements of the energy kernel to \Fin. Since this corresponds to minimization of energy, it is naturally related to capacity. 
\end{enumerate}
Both of these notions are methods of specifying a \emph{unique} solutions to $\Lap u = f$ in some way. The disparity between $R^F$ and $R^W$ is thus explained in terms of boundary conditions on \Lap as an unbounded self-adjoint operator on \HE in Remark~\ref{rem:wired-vs-free-as-boundary-conditions}.

\pgap

To compute effective resistance in an infinite network, we will consider three notions of subnetwork: free, wired, and trace. (Strictly speaking, these may not actually be subnetworks of the original graph; see Definitions~\ref{def:free-resistance}, \ref{def:wired-resistance}, and \ref{def:trace-resistance} for the precise details.) 
  Throughout this section, we use $H$ to denote a finite subnetwork of \Graph, $\verts[H]$ to denote its vertex set, and $H^F$, $H^W$, and $H^\trc$ to denote the free, wired, and trace networks associated to $\verts[H]$ (these terms are defined in other sections below). 

\begin{defn}\label{def:relative-resistance}
  If $H$ is a finite subnetwork of \Graph which contains $x$ and $y$, define $R_H(x,y)$ to be the \emph{resistance distance from $x$ to $y$ as computed within $H$}. In other words, compute $R_H(x,y)$ by any of the equivalent formulas of Theorem~\ref{thm:effective-resistance-metric}, but extremizing over only those functions whose support is contained in $H$.
\end{defn}

\begin{defn}\label{def:full-subnetwork}
  Let $\verts[H] \ci \verts$. Then the \emph{full subnetwork} on \verts[H] has all the edges of \Graph for which both endpoints lie in \verts[H], with the same conductances. That is, $\cond^H = \cond^G|_{\verts[H] \times \verts[H]}$.
\end{defn}

The notation $\{G_k\}_{k=1}^\iy$ always denotes an \emph{exhaustion} of the infinite network $(\Graph,\cond)$, as in Definition~\ref{def:exhaustion-of-G}. Since $x$ and $y$ are contained in all but finitely many $G_k$, we may always assume that $x,y \in G_k$, $\forall k$. Also, we assume in this section that the subnetworks $\{G_k\}$ are \emph{full}, in the sense of Definition~\ref{def:full-subnetwork}.
This may not be completely necessary, but it simplifies the discussion in a couple of places, avoids ambiguity, and causes no loss of generality.

\subsubsection{Free resistance}

\begin{defn}\label{def:free-resistance}
  For any subset $\verts[H] \ci \verts$, the \emph{free subnetwork} $H^F$ is just the full subnetwork with vertices $\verts[H]$. That is, all edges of \Graph with endpoints in \verts[H] are edges of $H^F$, with the same conductances. 
  Let $R_{H^F}(x,y)$ denote the effective resistance between $x$ and $y$ as computed in $H^F$, as in Definition~\ref{def:relative-resistance}. 
  The \emph{free resistance} between $x$ and $y$ is defined to be
  \linenopax
  \begin{align}\label{eqn:def:free-resistance}
    R^F(x,y) := \lim_{k \to \iy} R_{G_k^F}(x,y),
  \end{align}
  where $\{G_k\}$ is any exhaustion of \Graph.
\end{defn}
  \glossary{name={$R^F(x,y)$},description={free resistance metric},sort=RF,format=textbf}

\begin{remark}\label{rem:free-resistance}
  The name ``free'' comes from the fact that this formulation is free of any boundary conditions or considerations of the complement of $H$, in contrast to the wired and trace formulations of the next two subsections; see \cite[\S9]{Lyons:ProbOnTrees}.

  One can see that $R_{H^F}(x,y)$ has the drawback of ignoring the conductivity provided by all paths from $x$ to $y$ that pass through the complement of $H$. This provides some motivation for the wired and trace approaches below.
\end{remark}

\begin{defn}\label{def:L_xy}
  Fix $x, y \in \Graph$ and define the operator $L_{xy}$ on \HE by $L_{xy}v := v(x)-v(y)$. 
\end{defn}
We were led to \eqref{eqn:def:R(x,y)-sup} by considering the evaluation operators $L_x$ corresponding to $v_x$ in \cite[\S2.1]{DGG}; it is clear that \eqref{eqn:def:R(x,y)-sup} is equivalent to \eqref{eqn:def:R(x,y)-S} by considering the norm of $L_{xy} = L_x - L_y$.

\begin{remark}\label{rem:R(x,y) = |L{xy}|}
  Theorem~\ref{thm:free-resistance} is the free extension of Theorem~\ref{thm:effective-resistance-metric} to infinite networks; it shows that $R(x,y) = \|L_{xy}\|$ and that $R(x,o)$ is the best possible constant $k=k_x$ in \cite[Lemma~2.5]{DGG}.
\end{remark}

\begin{defn}\label{def:paths}
  A \emph{(finite) path} \cpath from $x \in \verts$ to $y \in \verts$ is a sequence of adjacent vertices $(x = x_0, x_1, x_2, \dots, x_n = y)$, i.e., $x_i \nbr x_{i-1}$ for $i=1,\dots,n$. The collection of paths from $x$ to $y$ is denoted $\Paths(x,y)$. 

  The \emph{characteristic function of a path} \cpath is the current $\charfn{\cpath}$ for which $\charfn{\cpath}(x_{i-1},x_i) = 1$, $i=1,\dots,n$, and $\charfn{\cpath}(x,y) = 0$ when $x$ and $y$ are not adjacent elements of \cpath. The notation $\curr = \sum_{\cpath \in \Paths(x,y)} \gx_\cpath \charfn{\cpath}$ indicates that \curr decomposes as a sum of currents supported on paths from $x$ to $y$. It may be that many (but not all) of the coefficients $\gx_\cpath$ are 0.
\end{defn}

\begin{theorem}\label{thm:free-resistance}
  For an infinite network \Graph, the free resistance $R^F(x,y)$ has the following equivalent formulations:
  \linenopax
  \begin{align}
    R^F(x,y)
    &= (v_x(x) - v_x(y)) - (v_y(x) - v_y(y)) \label{eqn:def:R^F(x,y)-Lap} \\
    &= \energy(v_x-v_y) \label{eqn:def:R^F(x,y)-energy} \\
    &= \min \{\diss(\curr) \suth \curr \in \Flo(x,y) \text{ and } \curr = \textstyle \sum_{\cpath \in \Paths(x,y)} \gx_\cpath \charfn{\cpath}\} \label{eqn:def:R^F(x,y)-diss} \\
    &= \tfrac1{\min\{\energy(u) \suth u(x)=1, u(y)=0, u \in \dom \energy\}} \label{eqn:def:R^F(x,y)-R} \\ 
    &= \inf\{\gk \geq 0 \suth |v(x)-v(y)|^2 \leq \gk \energy(v), v \in \dom \energy\} \label{eqn:def:R^F(x,y)-S} \\
    &= \sup\{|v(x)-v(y)|^2 \suth \energy(v) \leq 1, v \in \dom \energy\} \label{eqn:def:R^F(x,y)-sup} 
  \end{align}
  \begin{proof}
    To see that \eqref{eqn:def:R^F(x,y)-energy} is equivalent to \eqref{eqn:def:free-resistance}, fix any exhaustion of \Graph and note that 
    \linenopax
    \begin{align*}
      \energy(v_x-v_y) = \lim_{k \to \iy} \frac12 \sum_{s,t \in G_k} \cond_{st}((v_x-v_y)(s) - (v_x-v_y)(t))^2 = \lim_{k \to \iy} R_{G_k^F}(x,y),
    \end{align*}
    where the latter equality is from Theorem~\ref{thm:effective-resistance-metric}. Then for the equivalence of formulas \eqref{eqn:def:R^F(x,y)-Lap} and \eqref{eqn:def:R^F(x,y)-energy}, simply compute
    \linenopax
    \begin{align*}
      \energy(v_x-v_y)
      &= \la v_x-v_y, v_x-v_y \ra_\energy
       = \la v_x, v_x \ra_\energy - 2\la v_x, v_y \ra_\energy + \la v_y, v_y \ra_\energy
    \end{align*}
    and use the fact that $v_x$ is \bR-valued; cf.~\cite[Lemma~2.22]{DGG}.
    
    To see \eqref{eqn:def:R^F(x,y)-diss} is equivalent to \eqref{eqn:def:free-resistance}, fix any exhaustion of \Graph and define
    \linenopax
    \begin{align*}
      \Flo(x,y)\evald{H} := \{\curr \in \Flo(x,y) \suth \curr = \textstyle\sum\nolimits_{\cpath \ci H} \gx_\cpath \charfn{\cpath}\}.
    \end{align*}
    From \eqref{eqn:def:R(x,y)-diss}, it is clearly true for each $G_k$ that
    \linenopax
    \begin{align*}
      R_{G_k^F}(x,y) = \min \{\diss(\curr) \suth \curr \in \Flo(x,y) \text{ and } \curr = \textstyle \sum_{\cpath \ci G_k} \gx_\cpath \charfn{\cpath} \}.
    \end{align*}
    Since $\Flo(x,y)\evald{\Graph} = \bigcup_k \Flo(x,y)\evald{G_k}$, formula \eqref{eqn:def:R^F(x,y)-diss} follows. Note that \diss is a quadratic form on the closed convex set $\Flo(x,y)\evald{\Graph}$ and hence it attains its minimum.
    
    The equivalence of \eqref{eqn:def:R^F(x,y)-R} and \eqref{eqn:def:R^F(x,y)-sup} is \cite[Thm.~2.3.4]{Kig01}. 
    
    As for \eqref{eqn:def:R^F(x,y)-S} and \eqref{eqn:def:R^F(x,y)-sup}, they are both clearly equal to $\|L_{xy}\|$ (as described in Remark~\ref{rem:R(x,y) = |L{xy}|}) by the definition of operator norm; see \cite[\S5.3]{Rud87}, for example. To show that these are equivalent to $R^F$ as defined in \eqref{eqn:def:free-resistance}, define a subspace of \HE consisting of those voltages whose induced currents are supported in a finite subnetwork $H$ by    
    \linenopax
    \begin{equation}\label{eqn:def:HE-F-restr}
      \HE\evald[F]{H} = \{u \in \dom \energy \suth u(x)-u(y)=0 \text{ unless } x,y \in H\}. 
    \end{equation}
    This is a closed subspace, as it is the intersection of the kernels of a collection of continuous linear functionals $\|L_{st}\|$, and so we can let $Q_k$ be the projection to this subspace. Then it is clear that $Q_k \leq Q_{k+1}$ and that $\lim_{k \to \iy} \|u - Q_k u\|_\energy = 0$ for all $u \in \HE$, so
    \linenopax
    \begin{align}\label{eqn:defn-Q_k}
      R_{G_k^F}(x,y) 
      &= \| L_{xy} \|_{\HE|_{G_k} \to \bC }
       = \| L_{xy} Q_k\|,
    \end{align}
where the first equality follows from \eqref{eqn:def:R(x,y)-S} (recall that $G_k$ is finite) and therefore
    \linenopax
    \begin{align*}
      R^F(x,y) = \lim_{k \to \iy} R_{G_k^F}(x,y) 
      &= \lim_{k \to \iy} \|L_{xy} Q_k\| 
      = \left\| \lim_{k \to \iy} L_{xy} Q_k \right\| 
      = \|L_{xy}\|.
      \qedhere
    \end{align*}
  \end{proof}
\end{theorem}

  In view of the previous result, the free case corresponds to consideration of only those voltage functions whose induced current can be decomposed as a sum of currents supported on paths in \Graph. The wired case considered in the next section corresponds to considering all voltages functions whose induced current flow satisfies Kirchhoff's law in the form \eqref{eqn:Kirchhoff-nonhomog}; this is clear from comparison of \eqref{eqn:def:R^F(x,y)-diss} to \eqref{eqn:def:R^W(x,y)-diss}. See also Remark~\ref{rem:wired-vs-free-as-Kirchhoff-conditions+}.

  Formula \eqref{eqn:def:R^F(x,y)-Lap} turns out to be useful for explicit computations. Explicit formulas for the effective resistance metric on \bZd are obtained from \eqref{eqn:def:R^F(x,y)-Lap} in \cite[\S14.2]{OTERN}; compare to \cite[\S{V.2}]{Soardi94}. 

\begin{remark}\label{rem:martingales-in-Hilbert-space}
  In Theorem~\ref{thm:free-resistance}, the proofs that $R^F$ is given by \eqref{eqn:def:R^F(x,y)-diss} or \eqref{eqn:def:R^F(x,y)-S} stem from essentially the same underlying martingale argument. In a Hilbert space, a martingale is an increasing sequence of projections $\{Q_k\}$ with the martingale property $Q_k = Q_k Q_{k+1}$. Recall that conditional expectation is a projection. In this context, Doob's theorem \cite{Doob53} then states that if $\{f_k\} \ci \sH$ is such that $f_k = Q_k f_j$ for any $j \geq k$, then the following are equivalent:
\begin{enumerate}[(i)]
  \item there is a $f \in \sH$ such that $f_k = Q_k f$ for all $k$
  \item $\sup_k \|f_k\| < \iy$. 
\end{enumerate}
  The argument for \eqref{eqn:def:R^F(x,y)-diss} corresponds to projecting to subspaces of the Hilbert space of currents for which $\diss(\curr) < \iy$. 
  In \cite[\S9.1]{Lyons:ProbOnTrees}, the free resistance $R^F(x,y)$ is defined directly  via this approach (and similarly for $R^W(x,y)$). 
\end{remark}

The following result is also a special case of \cite[Thm.~2.3.4]{Kig01}.

\begin{prop}\label{thm:free-resistance-is-a-metric}
  $R^F(x,y)$ is a metric.
  \begin{proof}
    One has $R_{G_k^F}(x,z) \leq R_{G_k^F}(x,y) + R_{G_k^F}(y,z)$ for any $k$, so take the limit.
  \end{proof}
\end{prop}

\begin{cor}\label{thm:vx-is-Lipschitz}
  Any representative of $v \in \HE$, considered as a function on the metric space $(\Graph,R^F)$, is H\"{o}lder continuous with exponent $\frac12$ and constant $\|v\|_\energy$.
\end{cor}

The previous corollary is a restatement of \eqref{eqn:def:R^F(x,y)-S}, combined with the fact that $R^F(x,y)$ is always finite (which follows from connectedness of \Graph).
The Gaussian measure of Brownian motion is supported on the space of such functions \cite{Nelson64} and this is used in \cite{bdG}. Also, it is pointed out in \cite[Thm.~4.5]{Kig03} that 

\subsubsection{Wired resistance}

  \begin{figure}
    \centering
    \includegraphics{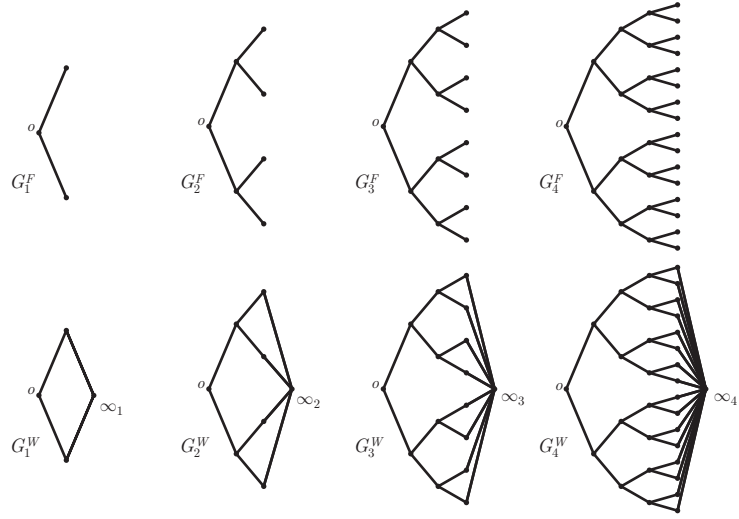}
    \caption{\captionsize Comparison of free and wired exhaustions for the example where \Graph is the infinite binary tree; see Definition~\ref{def:free-resistance} and Definition~\ref{def:wired-resistance} and also Example~\ref{exm:Binary-tree}. Here, the vertices of $G_k$ are all those which lie within $k$ edges (``steps'') of the origin. If the edges of \Graph all have conductance 1, then so do all the edges of each $G_k^F$ and $G_k^W$, except for edges incident upon $\iy_k = \iy_{G_k}$, which have conductance 2.}
    \label{fig:tree-exhaustions}
  \end{figure}

\begin{defn}\label{def:wired-resistance}
  Given a finite full subnetwork $H$ of \Graph, define the wired subnetwork $H^W$ by identifying all vertices in $\verts \less \verts[H]$ to a single, new vertex labeled \iy; see Figure\ref{fig:tree-exhaustions}. Thus, the vertex set of $H^W$ is $\verts[H] \cup \{\iy_H\}$, and the edge set of $H^W$ includes all the edges of $H$, with the same conductances. However, if $x \in \verts[H]$ has a neighbour $y \in \verts \less \verts[H]$, then $H^W$ also includes an edge from $x$ to \iy with conductance
  \linenopax
  \begin{align}\label{eqn:cond-to-iy}
    \cond_{x\iy_{\negsp[2]\scalebox{0.40}{$H$}}} := \sum_{y \nbr x, \, y \in H^\complm} \negsp[13]\cond_{xy}.
  \end{align}
  The identification of vertices in $G_k^\complm$ may result in parallel edges; then \eqref{eqn:cond-to-iy} corresponds to replacing these parallel edges by a single edge according to the usual formula for resistors in parallel.

  Let $R_{H^W}(x,y)$ denote the effective resistance between $x$ and $y$ as computed in $H^W$, as in Definition~\ref{def:relative-resistance}.
  The \emph{wired resistance} is then defined to be
  \linenopax
  \begin{align}\label{eqn:def:wired-resistance}
    R^W(x,y) := \lim_{k \to \iy} R_{G_k^W}(x,y),
  \end{align}
  where $\{G_k\}$ is any exhaustion of \Graph.
\end{defn}
  \glossary{name={$R^W(x,y)$},description={wired resistance metric},sort=RW,format=textbf}
  \glossary{name={$R_H(x,y)$},description={relative resistance as computed with respect to the subnetwork $H$},sort=RH}

\begin{remark}\label{rem:wired-resistance}
  The wired subnetwork is equivalently obtained by ``shorting together'' all vertices of $H^\complm$, and hence it follows from Rayleigh's monotonicity principle that $R^W(x,y) \leq R^F(x,y)$; cf. \cite[\S1.4]{DoSn84} or \cite[\S2.4]{Lyons:ProbOnTrees}.
\end{remark}

\begin{theorem}\label{thm:wired-resistance}
  The wired resistance may be computed by any of the following equivalent formulations:
  \linenopax
  \begin{align}
    R^W(x,y)
    &= \min_v \{v(x)-v(y) \suth \Lap v = \gd_x-\gd_y, v \in \dom\energy\}
    \label{eqn:def:R^W(x,y)-Lap} \\
    &= \min_v \{\energy(v) \suth \Lap v = \gd_x-\gd_y, v \in \dom\energy\} \label{eqn:def:R^W(x,y)-energy} \\
    &= \min_\curr \{\diss(\curr) \suth \curr \in \Flo(x,y), \diss(\curr) < \iy\} \label{eqn:def:R^W(x,y)-diss} \\
    &= 1/\min\{\energy(u) \suth u(x)=1, u(y)=0, u \in \Fin\} \label{eqn:def:R^W(x,y)-R} \\
    &= \inf\{\gk \geq 0 \suth |v(x)-v(y)|^2 \leq \gk \energy(v), v \in \Fin\} \label{eqn:def:R^W(x,y)-S} \\
    &= \sup\{|v(x)-v(y)|^2 \suth \energy(v) \leq 1, v \in \Fin\} \label{eqn:def:R^W(x,y)-sup} 
  \end{align}
  \begin{proof}
    Since \eqref{eqn:def:R^W(x,y)-S} and \eqref{eqn:def:R^W(x,y)-sup} are both clearly equivalent to the norm of $L_{xy}:\Fin \to \bC$ (where again $L_{xy}u - u(x)-u(y)$ as in Remark~\ref{rem:R(x,y) = |L{xy}|}), we begin by equating them to \eqref{eqn:def:wired-resistance}. From Definition~\ref{def:Fin}, we see that 
    \linenopax
    \begin{align}\label{eqn:def:HE-W-restr}
      \HE\evald[W]{H}
      := \{u \in \HE \suth \spt u \ci H\}
    \end{align}
    is a closed subspace of \HE. Let $Q_k$ be the projection to this subspace. Then it is clear that $Q_k \leq Q_{k+1}$ and that $\lim_{k \to \iy} \| P_{\Fin}u - Q_k u\|_\energy = 0$ for all $u \in \HE$. Each function $u$ on $H^W$ corresponds to a function $\tilde u$ on \Graph whose support is contained in $H$; simply define 
    \linenopax
    \begin{align*}
    \tilde u(x) = 
    \begin{cases} 
       u(x), & x \in H, \\ 
       u(\iy_H), & x \notin H.
    \end{cases}
    \end{align*}
    It is clear that this correspondence is bijective, and that
    \linenopax
    \begin{align*}
      R_{G_k^W}(x,y) 
      &= \| L_{xy} \|_{\HE|_{G_k}^W \to \bC }
       = \| L_{xy} Q_k\|,
    \end{align*}
    where the first equality follows from \eqref{eqn:def:R(x,y)-S} (recall that $G_k$ is finite) and therefore
    \linenopax
    \begin{align*}
      R^W(x,y) = \lim_{k \to \iy} R_{G_k^W}(x,y) 
      &= \lim_{k \to \iy} \|L_{xy} Q_k\| 
      = \|L_{xy} P_{\Fin}\|,
    \end{align*}
    which is equivalent to \eqref{eqn:def:R^W(x,y)-S}.
    
    To see \eqref{eqn:def:R^W(x,y)-Lap} is equivalent to \eqref{eqn:def:R^W(x,y)-energy}, note that the minimal energy solution to $\Lap u = \gd_x - \gd_y$ lies in \Fin, since any two solutions must differ by a harmonic function. Let $u$ be a solution to $\Lap u = \gd_x - \gd_y$ and define $f = \Pfin u$. Then $f \in \Fin$ and $\Lap f = \gd_x - \gd_y$ implies 
    \linenopax
    \begin{align}
      \|f\|_\energy^2 
      = \sum_{z \in \verts} f(z) \Lap f(z) 
      = \sum_{z \in \verts} f(z) (\gd_x-\gd_y)(z)
      = f(x) - f(y).
      \label{eqn:thm:wired-resistance:eqn1}
    \end{align}

    To see \eqref{eqn:def:R^W(x,y)-Lap}$\leq$ \eqref{eqn:def:R^W(x,y)-S}, let \gk be the optimal constant from \eqref{eqn:def:R^W(x,y)-S}. If $u \in \Fin$ is the unique solution to $\Lap u = \gd_x - \gd_y$, then 
    \linenopax
    \begin{align*}
      \gk = \sup_{u \in \Fin} \left\{\frac{|u(x)-u(y)|^2}{\energy(u)}\right\}
      \geq \frac{|u(x)-u(y)|^2}{\energy(u)} 
      = u(x)-u(y),               
    \end{align*}
    where the last equality follows from $\energy(u) = u(x) - u(y)$, by the same computation as in \eqref{eqn:thm:wired-resistance:eqn1}.
    For the reverse inequality, note that with $L_{xy}$ as just above,
    \linenopax
    \begin{align*}
      \frac{|u(x)-u(y)|^2}{\energy(u)}
      = \left| L_{xy} \left(\tfrac{u}{\energy(u)^{1/2}}\right)\right|^2
      = \left| \left\la v_x-v_y, \tfrac{u}{\energy(u)^{1/2}}\right\ra_\energy\right|^2,               
    \end{align*}
    for any $u \in \Fin$. Note that Lemma~\ref{thm:wired-resistance-is-a-metric} allows one to replace $v_x$ by $f_x = \Pfin v_x$, whence 
    \linenopax
    \begin{align*}
      \frac{|u(x)-u(y)|^2}{\energy(u)}
      &\leq \energy(f_x-f_y) \energy\left(\tfrac{u}{\energy(u)^{1/2}}\right) 
      = \energy(f_x-f_y) 
    \end{align*}
    by Cauchy-Schwarz. The infimum of the left-hand side over nonconstant functions $u \in \Fin$ gives the optimal \gk in \eqref{eqn:def:R^W(x,y)-S}, and thus shows that \eqref{eqn:def:R^W(x,y)-S} $\leq$ \eqref{eqn:def:R^W(x,y)-energy}. 

    To see \eqref{eqn:def:R^W(x,y)-energy} is equivalent to \eqref{eqn:def:R^W(x,y)-diss}, recall that \curr minimizes \diss over $\Flo(x,y)$ if and only if $\curr = \drp u$ for $u$ which minimizes \energy over $\{v \in \dom \energy \suth \Lap v = \gd_x - \gd_y\}$; see \cite[Thm.~3.26]{OTERN}, for example. Apply this to $\curr = \drp f$, where $f = \Pfin u$ is the minimal energy solution to $\Lap u = \gd_x - \gd_y$. 
    
    The equivalence of \eqref{eqn:def:R^W(x,y)-R} and \eqref{eqn:def:R^W(x,y)-sup} is directly parallel to the finite case and may also be obtained from \cite[Thm.~2.3.4]{Kig01}.
  \end{proof}
\end{theorem}
 
The proof of the next result follows from the finite case, exactly as in Theorem~\ref{thm:free-resistance-is-a-metric}.
\begin{theorem}\label{thm:wired-resistance-is-a-metric}
  $R^W(x,y)$ is a metric.
\end{theorem}

\begin{remark}[$R^F$ vs. $R^W$ explained in terms of boundary conditions on \Lap]
  \label{rem:wired-vs-free-as-boundary-conditions}
  Observe that both spaces  
  \linenopax
  \begin{align*} 
    \HE\evald[F]{H} = \{u \in \HE &\suth u(x)-u(y)=0 \text{ unless } x,y \in H\}
    \qq\text{and}\\
    \HE\evald[W]{H} &= \{u \in \HE \suth \spt u \ci H\}
  \end{align*}
  consist of functions which have no energy outside of $H$. The difference is that if the complement of $H$ consists of several connected components, then $u \in \HE|_H^F$ may take a different constant value on each one; this is not allowed for elements of $\HE|^W_H$. Therefore, $\HE|_H^F$ corresponds to Neumann boundary conditions and $\HE|_H^W$ corresponds to Dirichlet boundary conditions. That is, from the proofs of Theorem~\ref{thm:free-resistance} and Theorem~\ref{thm:wired-resistance}, we see
  \begin{enumerate}
    \item $R_{H^F}(x,y) = u(x)-u(y)$ where $u$ is the Neumann solution to $\Lap u = \gd_x - \gd_y$, and
    \item $R_{H^W}(x,y) = u(x)-u(y)$ where $u$ is the Dirichlet solution to $\Lap u = \gd_x - \gd_y$. 
  \end{enumerate}
  
\end{remark}

\begin{remark}\label{rem:wired-vs-free-as-Kirchhoff-conditions+}
  While the wired subnetwork takes into account the conductivity due to all paths from $x$ to $y$ (see Remark~\ref{rem:free-resistance}), it is overzealous in that it may also include paths from $x$ to $y$ that do not correspond to any path in \Graph (see Remark~\ref{rem:martingales-in-Hilbert-space}). On an infinite network, this leads to current ``pseudo-flows'' in which some of the current travels from $x$ to \iy, and then from \iy to $y$; see Example~\ref{exm:geometric-integers} and Theorem~\ref{thm:fx-as-prob}. 
\end{remark}

\begin{remark}[Comparison with resistance forms]
  \label{rem:comparison-to-resistance-forms}
  In \cite[Def.~2.8]{Kig03}, a \emph{resistance form} is defined as follows: let $X$ be a set and let \energy be a symmetric quadratic form on $\ell(X)$, the space of allfunctions on $X$, and let \sF denote the domain of \energy. Then $(\energy, \sF)$ is a resistance form iff:
  \begin{enumerate}[(RF1)]
    \item \sF is a linear subspace of $\ell(X)$ containing the constant functions and \energy is nonnegative on \sF with $\energy(u) = 0$ iff $u$ is constant.
    \item $\sF/\sim$ is a Hilbert with inner product \energy, where $\sim$ is the equivalence relation defined on \sF by $u \sim v$ iff $u-v$ is constant.
    \item For any finite subset $V \ci U$ and for any $v \in \ell(V)$, there is $u \in \sF$ such that $u \evald{V} = v$.
    \item For any $p,q \in X$, the number
      \begin{align}\label{eqn:form-resistance}
        R_{\energy,\sF}(p,q) := \sup \left\{ \tfrac{|u(p)-u(q)|^2}{\energy(u)} 
          \suth u \in \sF, \energy(u) > 0\right\}
      \end{align}
      is finite. Then $R_{\energy,\sF}$ is called the \emph{effective resistance associated to the form} $(\energy,\sF)$.
    \item If $u \in \sF$, then $\cj{u}$ defined by $\cj{u}(x) := \min\{1,\max\{0,u(x)\}\}$ (the unit normal contraction of $u$, in the language of Dirichlet forms) is also in \sF.
  \end{enumerate}
  Upon comparison of \eqref{eqn:def:R^F(x,y)-R}--\eqref{eqn:def:R^F(x,y)-S} to \eqref{eqn:def:R^W(x,y)-R}--\eqref{eqn:def:R^W(x,y)-sup}, one can see that $R^F$ is the effective resistance associated to the resistance form $(\energy,\HE)$, and that $R^W$ is the effective resistance associated to the resistance form $(\energy,\Fin)$. We are grateful to Jun Kigami for pointing this out to us. See also Remark~\ref{rem:Rharm-has-no-form}.
\end{remark}

\subsection{Harmonic resistance}
\label{sec:Harmonic-resistance}

\begin{defn}\label{def:harmonic-resistance}
  For an infinite network $(G,\cond)$ define the \emph{harmonic resistance} between $x$ and $y$ by
  \linenopax
  \begin{align}\label{eqn:def:harmonic-resistance}
    R^\hrm(x,y) := R^F(x,y) - R^W(x,y).
  \end{align}
\end{defn}

The next result is immediate upon comparing Theorem~\ref{thm:free-resistance} and Theorem~\ref{thm:wired-resistance}.

\begin{theorem}\label{thm:harmonic-resistance}
  With $h_x = \Phar v_x$ as in Remark~\ref{rem:3-repkernels}, the harmonic resistance is equal to
  \linenopax
  \begin{align}
    R^\hrm(x,y)
    &= (h_x(x) - h_x(y)) - (h_y(x) - h_y(y)) \label{eqn:def:R^ha(x,y)-Lap} \\
    &= \energy(h_x-h_y) \label{eqn:def:R^ha(x,y)-energy} \\
    &= \tfrac1{\min\{\energy(v) \suth v(x)=1, v(y)=0\}}  - \tfrac1{\min\{\energy(f) \suth f(x)=1, f(y)=0, f \in \Fin\}} \label{eqn:def:R^ha(x,y)-R} \\
    &= \inf\{\gk \geq 0 \suth |h(x)-h(y)|^2 \leq \gk \energy(h), h \in \Harm\} \label{eqn:def:R^ha(x,y)-S} \\
    &= \sup\{|h(x)-h(y)|^2 \suth \energy(h) \leq 1, h \in \Harm\} \label{eqn:def:R^ha(x,y)-sup} 
  \end{align}
\end{theorem}

  \version{}{\marginpar{Idea for precise definition: $g_o$ is constant on this portion.}}
\begin{remark}\label{rem:harmonic-resistance-not-a-metric}\label{rem:Rharm-has-no-form}
  Note that $R^\hrm$ is not the effective resistance associated to a resistance form, as in Remark~\ref{rem:comparison-to-resistance-forms}, since (RF5) may fail. If $R^\hrm$ \emph{were} the effective resistance associated to a resistance form, then \cite[Prop.~2.10]{Kig03} would imply that $R^\hrm(x,y)$ is a metric, but this can be seen to be false by considering basic examples. See Example~\ref{exm:geometric-integers}, e.g. The same remarks also apply to the  boundary resistance $R^{\bdy}(x,y)$, discussed just below. 
\end{remark}

\begin{defn}\label{def:boundary-resistance}
  For an infinite network $(G,\cond)$ define the \emph{boundary resistance} between $x$ and $y$ by
  \linenopax
  \begin{align}\label{eqn:def:boundary-resistance}
    R^{\bdy}(x,y) := \frac1{R^W(x,y)^{-1} - R^F(x,y)^{-1}}.
  \end{align}
\end{defn}

  Intuitively, some portion of the wired/minimal current from $x$ to $y$ passes through infinity; the quantity $R^{\bdy}(x,y)$ gives the voltage drop ``across infinity''; see Remark~\ref{rem:shortcut-through-infinity}. From this perspective, infinity is ``connected in parallel''. The boundary $\bd G$ in \cite{bdG} is a more rigorous definition of the set at infinity.

\begin{theorem}\label{thm:boundary-resistance}
  The boundary resistance is equal to
  \linenopax
  \begin{align}
    R^\bdy(x,y) = 
    \frac{R^W(x,y)R^F(x,y)}{R^{\hrm}(x,y)}.
    \label{eqn:def:R^hrm(x,y)-R}
  \end{align}
  In particular, the resistance across the boundary is infinite if $\Harm = 0$. 
  \begin{proof}
    From \eqref{eqn:def:harmonic-resistance} one has $R^F(x,y) = 1/(R^W(x,y)^{-1} - R^\bdy(x,y)^{-1})$, which gives
    \begin{align*}
      \frac1{\energy(v_x-v_y)}
        &=  \frac1{\energy(f_x-f_y)} - \frac1{R^\bdy(x,y)}
    \end{align*}
    by Theorem~\ref{thm:free-resistance} and Theorem~\ref{thm:wired-resistance}, and hence
    \begin{align*}
      \frac1{R^\bdy(x,y)}
        &= \frac1{\energy(f_x-f_y)} - \frac1{\energy(v_x-v_y)}.
    \end{align*}
    Now solving for $R^\bdy$ gives
    \begin{align}
      R^\bdy(x,y) 
      = \frac{\energy(f_x-f_y)\energy(v_x-v_y)}{\energy(h_x-h_y)},
    \end{align}
    and the conclusion follows from \eqref{eqn:def:R^F(x,y)-energy}, \eqref{eqn:def:R^W(x,y)-energy}, and \eqref{eqn:def:R^ha(x,y)-energy}.
  \end{proof}
\end{theorem}


\section{Trace resistance}
\label{sec:trace-resistance}

The third type of subnetwork takes into account the connectivity of the complement of the subnetwork, but does not add anything extra.  The name ``trace'' is due to the fact that this approach comes by considering the trace of the Dirichlet form \energy to a subnetwork; see \cite{FOT94}. Several of the ideas in this section were explored previously in \cite{Kig01,Kig03,Kig09,Metz}.

The discussion of the trace resistance and trace subnetworks requires some definitions relating the transition operator (i.e. Markov chain) \Prob to the probability measure $\prob^{(\cond)}$ on the space of (infinite) paths in \Graph which start at a fixed vertex $a$. Such a path is a sequence of vertices $\{x_n\}_{n=0}^\iy$, where $x_0=a$ and $x_n \nbr x_{n+1}$ for all $n$. 
  
\begin{defn}\label{def:Gamma(a,b)}
  Let $\Paths(a)$ be the space of all paths \cpath beginning at the vertex $a \in \verts$, and let $\Paths(a,b) \ci \Paths(a)$ be the subset of paths that reach $b$, and that do so \emph{before} returning to $a$:
  \linenopax
  \begin{align}\label{eqn:def:Paths(a,b)}
    \Paths(a,b) := \{\cpath \in \Paths(a) \suth b=x_n \text{ for some $n$, with } x_k \neq a \text{ for } 1 \leq k \leq n\}.
  \end{align}
\end{defn}

\begin{defn}\label{def:prob^cond}
   The space $\Paths(a)$ carries a natural probability measure $\prob^{(\cond)}$ defined by 
  \linenopax
  \begin{align}\label{eqn:prob^cond}
    \prob^{(\cond)}(\cpath) := \prod_{x_i \in \cpath} p(x_{i-1},x_i),
  \end{align}
  where $p(x,y)$ is as in \eqref{eqn:def:Prob-trans-oper}.
  The construction of $\prob^{(\cond)}$ comes by applying Kolmogorov consistency to the natural cylinder-set Borel topology that makes $\Paths(a)$ into a compact Hausdorff space. 
\end{defn}

\begin{defn}\label{def:hitting-time}
  Let $X_m$ be a random variable which denotes the (vertex) location of the random walker at time $m$. Then let $\gt_x$ be the \emph{hitting time} of $x$, that is, the random variable which is the expected time at which the walker first reaches $x$:
  \begin{align}\label{eqn:def:hitting-time}
    \gt_x := \min\{m \geq 0 \suth X_m = x\}.
  \end{align}
  More generally, $\gt_H$ is the time at which the walker first reaches the subnetwork $H$. For a walk started in $H$, this gives $\gt_H = 0$. 
\end{defn}

\subsection{The trace subnetwork}

It is well-known that networks $\{(\Graph,\cond)\}$ are in bijective correspondence with reversible Markov processes $\{\Prob\}$; this is immediate from the \emph{detailed balance equations} which follow from the symmetry of the conductance:
\begin{align*}
  \cond(x) p(x,y) = \cond_{xy} = \cond_{yx} = \cond(y) p(y,x).
\end{align*}
It follows from $\Lap = \cond(\one-\Prob)$ that networks are thus in bijective correspondence with \emph{Laplacians}, if one defines a Laplacian as in \eqref{eqn:def:laplacian}. That is, a Laplacian is a symmetric linear operator which is nonnegative definite, has kernel consisting of the constant functions, and satisfies $(\Lap\gd_x)(y) \leq 0$ for $x \neq y$. In other words, every row (and column) of $\trc(\Lap,H)$ sums to 0. (This is the negative of the definition of a \emph{Laplacian} as in \cite{Kig01} and \cite{CdV98}.) 
In this section, we exploit the bijection between Laplacians and networks to define the trace subnetwork. For $\verts[H] \ci \Graph$, the idea is as follows:
\begin{align*}
  G \;\longleftrightarrow\; \Lap 
  \;\limmode[\text{take the trace to \verts[H]}]\; \trc(\Lap, \verts[H]) 
  \;\longleftrightarrow\; H^\trc.
\end{align*}

\begin{defn}\label{def:trace-resistance}
  The \emph{trace} of \Graph to $\verts[H]$ is the network whose edge data is defined by the trace of \Lap to \verts[H], which is computed as the Schur complement of the Laplacian of $H$ with respect to \Graph. More precisely, write the Laplacian of \Graph as a matrix in block form, with the rows and columns indexed by vertices, and order the vertices so that those of $H$ appear first:
  \linenopax
  \begin{align}\label{eqn:Lap-block-decomp}
    \Lap =
      \begin{array}{l} \scalebox{0.70}{$H$} \\ \scalebox{0.70}{$H^\complm$} \end{array}
      \negsp[12]
      \left[\begin{array}{ll} A & B^T \\ B & D \end{array}\right],
  \end{align}
  where $B^T$ is the transpose of $B$. If $\ell(G) :=\{f:\verts \to \bR\}$, the corresponding mappings are
  \linenopax
  \begin{align}
    A:&\ell(H) \to \ell(H) & B^T&:\ell(H^\complm) \to \ell(H) \notag \\
    B:&\ell(H) \to \ell(H^\complm) & D&:\ell(H^\complm) \to \ell(H^\complm).
    \label{eqn:action-of-Lap-quadrants}
  \end{align}
  It turns out that the Schur complement
  \linenopax
  \begin{equation}\label{eqn:def:Schur-complement}
    \trc(\Lap, \verts[H]) := A - B^T D^{-1} B
  \end{equation}
  is the Laplacian of a subnetwork with vertex set $\verts[H]$; cf.~\cite[\S2.1]{Kig01} and Remark~\ref{rem:trace-valid-for-subsets}.\footnote{It will be clear from \eqref{eqn:Schur-complement-as-sum} that $D^{-1}$ always exists in this context, and hence \eqref{eqn:def:Schur-complement} is always well-defined. Furthermore, the existence of the trace is given in \cite[Prop.~2.10]{Kig03}; it is known from \cite[Lem.~2.1.5]{Kig01} that $D$ is invertible and negative semidefinite.} A formula for the conductances (and hence the adjacencies) of the trace is given in Theorem~\ref{thm:edges-of-H^S}. Denote this new subnetwork by $H^\trc$.

  If $\verts[H] \ci \verts$ is finite, then for $x,y \in H$, the trace of the resistance on $H$ is denoted $R_{H^\trc}(x,y)$, and defined as in Definition~\ref{def:effective-resistance}. The \emph{trace resistance} is then defined to be
  \linenopax
  \begin{align}\label{eqn:trace-resistance}
    R^\trc(x,y) := \lim_{k \to \iy} R_{G_k^\trc}(x,y),
  \end{align}
  where $\{G_k\}$ is any exhaustion of \Graph.
\end{defn}  

\begin{defn}\label{def:paths-through-the-complement}
    If $a,b \in \bd H$, then we write 
  \linenopax
  \begin{align}\label{eqn:def:Paths(a,b)-outside-H}
    \Paths(a,b)\evald{H^\complm} 
    := \{\cpath \in \Paths(a,b) \suth x_i \in H^\complm, 0 < i < \gt_b\},
  \end{align}
  for the set of paths from $a$ to $b$ that do not pass through any vertex in $\verts[H]$. 
  Note that if $x,y \in \bd H$ are adjacent, then any path of the form $\cpath = (x,y,\dots)$ is trivially in $\Paths(a,b)\evald{H^\complm}$. 
\end{defn}
  
\begin{defn}\label{def:prob[a->b]}
  Let $\prob[a \to b]$ denote the probability that a random walk started at $a$ will reach $b$ before returning to $a$. That is,
  \glossary{name={$\prob[a \to b]$},description={probability of the random walk started from $a$ reaching $b$},sort=P,format=textbf}
  \linenopax
  \begin{align}\label{eqn:def:prob-a-to-b}
    \prob[a \to b] := \prob^{(\cond)}(\Paths(a,b)).
  \end{align}
  Note that this is equivalent to 
  \linenopax
  \begin{align}\label{eqn:def:prob-a-to-b-as-hitting}
    \prob[a \to b] = \prob_a[\gt_b < \gt_a^+] := \prob[\gt_b < \gt_a \,|\, x_0=a],
  \end{align}
  where $\gt_a^+$ is the \emph{hitting time of $a$ after leaving $a$}, i.e., $\gt_a^+ := \min\{m \geq 1 \suth X_m = a\}$; compare to Definition~\ref{def:hitting-time}. More generally, one also has $\gt_H^+ := \min\{m \geq 1 \suth X_m \in H\}$.
  
  If $a,b \in \bd H$, then we write 
  \linenopax
  \begin{align}\label{eqn:def:prob-a-to-b-outside-H}
    \prob[a \to b]\evald{H^\complm} 
    := \prob^{(\cond)}\left(\Paths(a,b)\evald{H^\complm}\right),
  \end{align}
  \glossary{name={$\complement$},description={$H^\complm$ is the complement of $H$ in $G$},sort=C}
  that is, the probability that a random walk started at $a$ will reach $b$ via a path for which $X_m \notin H$ for $m=1,2,\dots,\gt_b-1$.
\end{defn}

\begin{remark}[More probabilistic notation]\label{rem:prob-notation}\label{rem:notation-restriction-to-Gk}
  The formulation in \eqref{eqn:def:prob-a-to-b-outside-H} is conditioning $\prob^{(\cond)}(\Paths(a,b))$ on avoiding $H$; the notation is intended to evoke something like ``$\prob[a \to b \,|\, \cpath \ci H^\complm]$''. However, this is not correct because $a,b \in H$ and \cpath may pass through $H$ after $\gt_b$.

  In Theorem~\ref{thm:edges-of-H^S}, we use the following common notation as in \cite{Spitzer} or \cite{Woess00}, for example. All notations are for the random walk started at $x$.
  \begin{align*}
    \Prob^n(x,y) = p^{(n)}(x,y) =\prob_x[X_n = y] &\q\text{probability that the walk is at $y$ after $n$ steps} \\
    G(x,y) = {\textstyle\sum}_{n=0}^\iy p^{(n)}(x,y) &\q\text{expected number of visits to $y$} \\
    f^{(n)}(x,y) = \prob_x[\gt_y = n] &\q\text{probability that the walk first reaches $y$ on the \nth step} \\
    F(x,y) = {\textstyle\sum}_{n=0}^\iy f^{(n)}(x,y) &\q\text{probability that the walk ever reaches $y$}
  \end{align*}
  Note that if the walk is killed when it reaches $y$, then $p^{(n)}(x,y) = f^{(n)}(x,y)$ because the first time it reaches $y$ is the only time it reaches $y$. Therefore, when the walk is conditioned to end upon reaching a set $S$, one has $G(x,y) = F(x,y)$ for all $y \in S$.
\end{remark}

\begin{theorem}\label{thm:edges-of-H^S}
  For $\verts[H] \ci \verts$, the conductances in the trace subnetwork $H^\trc$ are given by
  \linenopax
  \begin{align}\label{eqn:conductances-in-HS}
    \cond_{xy}^\trc = \cond_{xy} + \cond(x) \prob[x \to y]\evald{H^\complm}.
  \end{align}
  Consequently, the transition probabilities in the trace subnetwork are given by
  \linenopax
  \begin{align}\label{eqn:transition-probabilities-in-HS}
    p^\trc(x,y) = p(x,y) + \prob[x \to y]\evald{H^\complm}.
  \end{align}
\begin{proof}
  Using subscripts to indicate the block decomposition corresponding to $H$ and $H^\complm$ as in \eqref{eqn:Lap-block-decomp}, the Laplacian may be written as
  \linenopax
  \begin{align*}
    \Lap =
    \left[\begin{array}{cc} \cond_A(\one - \Prob_A) & -\cond_A\Prob_{B^T} \\ -\cond_D\Prob_{B} & \cond_D(\one - \Prob_{D}) \end{array}\right],
    \qq \text{for} \qq
    \cond =
      \begin{array}{l} \scalebox{0.70}{$H$} \\ \scalebox{0.70}{$H^\complm$} \end{array}
      \negsp[12]
      \left[\begin{array}{ll} \cond_A &  \\  & \cond_D \end{array}\right].
  \end{align*}
  Then the Schur complement is
  \linenopax
  \begin{align}
    \trc(\Lap,H)
    &= \cond_A - \cond_A \Prob_A - \cond_A \Prob_{B^T} (\id - \Prob_D)^{-1} \cond_{D}^{-1} \cond_{D} \Prob_B \notag \\
    &= \cond_A - \cond_A\left(\Prob_A + \Prob_{B^T} \left(\sum_{n=0}^\iy \Prob_D^n\right) \Prob_B \right) \notag \\
    &= \cond_A(\id - \Prob_\xX).
    \label{eqn:Schur-complement-as-sum}
  \end{align}
  Note that $\Prob_D$ is substochastic, and hence the RW has positive probability of hitting $\bd G_k$, whose vertices act as absorbing states. This means that the expected number of visits to any vertex in $H^\complm$ is finite and hence the matrix $\Prob_\xX$ has finite entries. 
  
  Meanwhile, using $\Prob_A(x,y)$ to denote the \nth[(x,y)] entry of the matrix $\Prob_A$, and $\gt_H^+$ as in Definition~\ref{def:prob[a->b]}, we have
  \linenopax
  \begin{align}
    \prob[x \to y]\evald{H^\complm}
    &= \prob^{(\cond)}\left(\Paths(x,y)\evald{H^\complm}\right) \notag \\
    &= \prob^{(\cond)}\left(\bigcup_{k=1}^\iy \{\cpath \in \Paths(x,y)\evald{H^\complm} \suth \gt_H^+ = k \}\right) \notag \\
    &= \prob^{(\cond)}\left(\{\cpath \in \Paths(x,y)\evald{H^\complm} \suth \gt_H^+ = 1 \}\right)
      + \sum_{k=2}^\iy \prob^{(\cond)}\left(\{\cpath \in \Paths(x,y)\evald{H^\complm} \suth \gt_H^+ = k \}\right) \notag \\
    &= \Prob_A(x,y) +\sum_{n=0}^\iy \sum_{s,t} \Prob_{B^T}(x,s) \Prob_D^n(s,t) \Prob_B(t,y) \label{eqn:prob-a-to-b-avoiding-H-computation}\\ 
    &= \Prob_\xX(x,y). \notag
  \end{align}
  To justify \eqref{eqn:prob-a-to-b-avoiding-H-computation}, note that by \eqref{eqn:action-of-Lap-quadrants}, $\Prob_D^n$ corresponds to steps taken in $H^\complm$. Therefore,
  \begin{align*}
    \left(\Prob_{B^T} \left(\sum_{n=0}^\iy \Prob_D^n\right) \Prob_B\right)(x,y)
    = \Prob_{B^T}\Prob_B(x,y) + \Prob_{B^T}\Prob_D \Prob_B(x,y) 
     + \Prob_{B^T}\Prob_D^2 \Prob_B(x,y) + \dots
  \end{align*}
  is the probability of the random walk taking a path that steps from $x \in H$ to $H^\complm$, meanders through $H^\complm$ for any finite number of steps, and finally steps to $y \in H$. Since $y \notin H^\complm$, 
  \begin{align*}
    \Prob_{B^T} \Prob_D^k \Prob_B(x,y) 
    = \prob_x[X_{k+2} = y] 
    = \prob_x[\gt_y = k+2] ,
  \end{align*}
because the walk can only reach $y$ on the last step, as in Remark~\ref{rem:prob-notation}. It follows by classical theory (see \cite{Spitzer}, for example) that the sum in \eqref{eqn:prob-a-to-b-avoiding-H-computation} is a probability (as opposed to an expectation, etc.) and justifies the probabilistic notation $\Prob_\xX$ in \eqref{eqn:Schur-complement-as-sum}.
  Note that $\Prob_A(x,y)$ corresponds to the one-step path from $x$ to $y$, which is trivially in $\Paths(x,y)\evald{H^\complm}$ by \eqref{eqn:def:Paths(a,b)-outside-H}.
  Since $\Prob_A(x,y) = p(x,y) = \cond_{xy}/\cond(x)$, the desired conclusion \eqref{eqn:conductances-in-HS} follows from combining \eqref{eqn:Schur-complement-as-sum}, \eqref{eqn:prob-a-to-b-avoiding-H-computation}, and \eqref{eqn:def:prob-a-to-b-outside-H}. Of course, \eqref{eqn:transition-probabilities-in-HS} follows immediately by dividing through by $\cond(x)$.
\end{proof}
\end{theorem}
 
  The authors are grateful to Jun Kigami for helpful conversations and suggestions for the proof of Theorem~\ref{thm:edges-of-H^S}. 

\begin{remark}[The trace construction is valid for general subsets of vertices]
  \label{rem:trace-valid-for-subsets}\label{rem:edges-of-HS}
  While Definition~\ref{def:trace-resistance} applies to a (connected) subnetwork of \Graph, it is essential to note that Theorem~\ref{thm:edges-of-H^S} applies to arbitrary subsets \verts[H] of \verts. 
  
  It is clear from \eqref{eqn:conductances-in-HS} that the edge sets of $\inn H$ and $\inn H^\trc$ are identical, but the conductance between two vertices $x,y \in \bd H^\trc$ is greater iff there is a path from $x$ to $y$ that does not pass through $H$. Indeed, if there is a path from $x$ to $y$ which lies entirely in $H^\complm$ except for the endpoints, then $x$ and $y$ will be adjacent in $H^\trc$, even if they were not adjacent in $H$.
\end{remark}

\begin{remark}[Network reduction, and resistance distance via Schur complement]
  \label{rem:Schur-does-network-reduction}
  \label{rem:Resistance-distance-via-Schur-complement} 
  A theorem of Epifanov states that every finite planar network with vertices $x,y$ can be reduced to a single equivalent conductor via the use of three simple transformations: parallel, series, and $\nabla$-\textsf{Y}; cf. \cite{Epifanov66,Truemper89} as well as \cite[\S2.3]{Lyons:ProbOnTrees} and \cite[\S7.4]{CdV98}. More precisely,
  \begin{enumerate}[(i)]
    \item Parallel. Two conductors $\cond_{xy}^{(1)}$ and $\cond_{xy}^{(2)}$ connected in parallel can be replaced by a single conductor $\cond_{xy} = \cond_{xy}^{(1)} + \cond_{xy}^{(2)}$.
    \item Series. If $z$ has only the neighbours $x$ and $y$, then $z$ may be removed from the network and the edges $\cond_{xz}$ and $\cond_{yz}$ should be replaced by a single edge $\cond_{xy} = (\cond_{xz}^{-1}+\cond_{yz}^{-1})^{-1}$.
    \item $\nabla$-\textsf{Y}. Let $t$ be a vertex whose only neighbours are $x,y,z$. Then this ``\textsf{Y}'' may be replaced by a triangle (``$\nabla$'') which does not include $t$, with conductances
      \linenopax
      \begin{align*}
        \cond_{xy} = \frac{\cond_{xt}\cond_{ty}}{\cond(t)}, \;
        \cond_{yz} = \frac{\cond_{yt}\cond_{tz}}{\cond(t)}, \;
        \cond_{xz} = \frac{\cond_{xt}\cond_{tz}}{\cond(t)}.
      \end{align*}
      This transformation may also be inverted, to replace a $\nabla$ with a \textsf{Y} and introduce a new vertex. 
  \end{enumerate}
  It is a fun exercise to obtain the series and $\nabla$-\textsf{Y} formulas by applying the Schur complement technique to remove a single vertex of degree 2 or 3 from a network. Indeed, these are both special cases of the following: let $t$ be a vertex of degree $n$, and let $H$ be the (star-shaped) subnetwork consisting only of $t$ and its neighbours. If we write the Laplacian for just this subnetwork with the \nth[t] row \& column last, then
  \linenopax
  \begin{align*}
    \Lap|_H =
    \left[\begin{array}{rrrc}
      \cond_{x_1t} & \dots & 0 & -\cond_{x_1t} \\
      \vdots & \ddots & \vdots & \vdots \\
      0 & \dots & \cond_{x_nt} & -\cond_{x_nt} \\
      -\cond_{x_1t} & \dots & -\cond_{x_nt} & \cond(t)
    \end{array}\right]
  \end{align*}
  and the Schur complement is
  \linenopax
  \begin{align*}
    \trc(\Lap|_H,H\less\{t\})
    =
    \left[\begin{array}{rrrc}
      \cond_{x_1t} & \dots & 0 \\
      \vdots & \ddots & \vdots \\
      0 & \dots & \cond_{x_nt} \\
    \end{array}\right]
    - \frac{1}{\cond(t)}
    \left[\begin{array}{rrrc}
      \cond_{x_1t} \\
      \vdots \\
      \cond_{x_nt} 
    \end{array}\right]
    \left[\begin{array}{rrrc}
      \cond_{x_1t} & \dots & \cond_{x_nt} 
    \end{array}\right],
  \end{align*}
  whence the new conductance from $x_i$ to $x_j$ is given by $\cond_{x_it}\cond_{tx_j}/\cond(t)$. It is interesting to note that the operator being subtracted corresponds to the projection to the rank-one subspace spanned by the probabilities of leaving $t$:
  \linenopax
  \begin{align*}
    \frac{1}{\cond(t)}
    \left[\begin{array}{r}
      \cond_{x_1t} \\
      \vdots \\
      \cond_{x_nt} 
    \end{array}\right]
    \left[\begin{array}{rrr}
      \cond_{x_1t} & \dots & \cond_{x_nt} 
    \end{array}\right]
    = \cond(t) |v\ra \la v|,
  \end{align*}
  using Dirac's ket-bra notation for the projection to a rank-1 subspace spanned by $v$ where
  \linenopax
  \begin{align*}
    v = \left[\begin{array}{rrr}
      p(t,x_1)  & \dots & p(t,x_n)  
    \end{array}\right].
  \end{align*}
In fact, $|v\ra \la v | = \Prob_\xX$, in the notation of \eqref{eqn:Schur-complement-as-sum}. In general, the trace construction (Schur complement) has the effect of probabilistically projecting away the complement of the subnetwork.

   Remark~\ref{rem:Resistance-distance-via-network-reduction} describes how the effective resistance can be interpreted as the correct resistance for a single edge which replaces a subnetwork; see Figure~\ref{fig:network-reduction}. The following corollary of Theorem~\ref{thm:edges-of-H^S} formalizes this interpretation by exploiting the fact that the Schur complement construction is viable for arbitrary subsets of vertices; see Remark~\ref{rem:trace-valid-for-subsets}. In this case, one takes the trace of the (typically disconnected) subset $\{x,y\} \ci \verts$; note that \raisebox{0.8mm}{\scalebox{0.50}{$\left[\begin{array}{rr} 1 & -1 \\ -1 & 1 \end{array}\right]$}} 
   is the Laplacian of the trivial 2-vertex network when the edge between them has unit conductance. The following result is \cite[Thm.~2.14]{Kig03}.
\end{remark}

\begin{lemma}
  \label{thm:trace-invariance-of-resistance}
  Let $H_2 \ci H_1$ be finite subnetworks of \Graph. Then for $a,b \in \verts[H_2]$, one has $R_{H_1^\trc}(a,b) = R_{H_2^\trc}(a,b)$.
\end{lemma}

\begin{cor}\label{thm:resistance-as-2-vert-network}
  Let $\verts[H] = \{x,y\}$ be any two vertices of \Graph. Then the trace resistance is
  \linenopax
  \begin{equation}\label{eqn:resistance-as-2-vert-network}
    \trc(\Lap,H)
    = \frac1{R^\trc(x,y)} \left[\begin{array}{rr} 1 & -1 \\ -1 & 1 \end{array}\right],
  \end{equation}
  where $\trc(\Lap,H) = A - B^T D^{-1} B$ as in \eqref{eqn:def:Schur-complement}.
  \begin{proof}
    Take $H=\{x,y\}$ in Theorem~\ref{thm:edges-of-H^S}. As discussed in Remark~\ref{rem:trace-valid-for-subsets}, it is not necessary to have $x \nbr y$. For any exhaustion $\{G_k\}_{k=1}^\iy$ with $H = G_1$, Corollary~\ref{thm:trace-invariance-of-resistance} gives 
    \begin{align*}
      R^\trc(x,y) = \lim_{k \to \iy} R_{G_k^\trc}(x,y)
       = R_{G_k^\trc}(x,y)
       = R_{H^\trc}(x,y).
    \end{align*}
    Note that in this case, $\left(\Prob_{B^T} \sum_n \Prob_D^n \Prob_B\right)(x,y)$ corresponds all paths from $x$ to $y$ that consist of more than one step:
  \linenopax
  \begin{align}
    \prob[x \to y]\evald{H^\complm}
    &= \Prob_A(x,y) + \left(\Prob_{B^T} \sum_{n=0}^\iy \Prob_D^n \Prob_B\right)(x,y) 
        = p(x,y) 
    + \sum_{|\cpath| \geq 2} \prob(\cpath),
    \label{eqn:sum-over-paths-decomposition}
  \end{align}
  and that the matrix in \eqref{eqn:resistance-as-2-vert-network} is the Laplacian of a network consisting of two vertices and one edge of conductance 1; see Figure~\ref{fig:network-reduction}.
  \end{proof}
\end{cor}

\version{}{\marginpar{Alternative idea}Use Kigami's approach (\energy, Theorem~\ref{thm:resistance-as-2-vert-network}) and take limits.}

\begin{cor}\label{thm:resistance-as-path-integral}
  The trace resistance $R^\trc(x,y)$ is given by
  \linenopax
  \begin{equation}\label{eqn:resistance-as-path-integral}
    R^\trc(x,y) = \frac1{\cond(x) \prob[x \to y]}
  \end{equation}
  \begin{proof}
    Again, take $\verts[H] = \{x,y\}$.  Then    
    \linenopax
    \begin{align*}
      R^\trc(x,y)^{-1}
      = \cond_{xy}^{H^\trc} 
      &= \cond_{xy} + \cond(x) \prob[x \to y]\evald{H^\complm} \\
      &= \cond(x) \left(p(x,y) + \prob[x \to y]\evald{H^\complm}\right) \\
      &= \cond(x) \prob[x \to y],
    \end{align*}
    where Corollary~\ref{thm:resistance-as-2-vert-network} gives the first equality and Theorem~\ref{thm:edges-of-H^S} gives the second.
  \end{proof}
\end{cor}
   
\begin{remark}[Effective resistance as ``path integral'']
  \label{rem:resistance-as-path-integral}
  Corollary~\ref{thm:resistance-as-path-integral} may also be obtained by the more elegant (and much shorter) approach of \cite[\S2.2]{Lyons:ProbOnTrees}, where it is stated as follows: the mean number of times a random walk visits $a$ before reaching $b$ is $\prob[a \to b]^{-1} = \cond(a) R(a,b)$. We give the present proof to highlight and explain the underlying role of the Schur complement with respect to network reduction; see Remarks~\ref{rem:trace-valid-for-subsets}--\ref{rem:Resistance-distance-via-Schur-complement}.
  A key point of the present approach is to emphasize the expression of effective resistance $R(a,b)$ in terms of a \emph{sum over all possible paths from $a$ to $b$}.   By Remark~\ref{rem:wired-vs-free-as-Kirchhoff-conditions+}, it is apparent that this ``path-integral'' interpretation makes $R^\trc$ much more closely related to $R^F$ than to $R^W$, as seen by the following result, which also follows immediately from \cite[Thm.~2.14]{Kig03}.
\end{remark}

\begin{cor}\label{thm:trace-resistance-is-free-resistance}
  On any transient network, $R^\trc(a,b) = R^F(a,b)$.
  \begin{proof}
    By Corollary~\ref{thm:trace-invariance-of-resistance}, it is clear that $R_{G_k^\trc}(a,b) = R_{G_{k+1}^\trc}(a,b)$ for all $k$. Meanwhile, any path from $a$ to $b$ will lie in $G_k$ for sufficiently large $k$, so it is clear by Theorem~\ref{thm:resistance-as-path-integral}, the sequence $\{R^F_{G_k}(a,b)\}_{k=0}^\iy$ is monotonically decreasing with limit $R^F(a,b) = R^\trc(a,b)$.
  \end{proof}
\end{cor}

\begin{remark}
  \version{}{\marginpar{Work out, and insert here, the corresponding formula for $R^W$}}
  Writing $[x \to y \,|\, \cpath \ci H]$ to indicate a restriction to paths from $x$ to $y$ that lie entirely in $H$, as in Remark~\ref{rem:notation-restriction-to-Gk}, one has
  \linenopax
  \begin{align*}
    R_{G_k^\trc}(x,y) 
    &= \frac1{\cond(x) \left(\prob[x \to y \,|\, \cpath \ci G_k] + \prob[x \to y  \,|\, \cpath \nsubseteq G_k]\right)}  \\
    &\leq \frac1{\cond(x) \prob[x \to y \,|\, \cpath \ci G_k]} 
    = R_{G_k^F}(x,y).
  \end{align*}
  Essentially, Corollary~\ref{thm:trace-invariance-of-resistance} is an expression of the first equality and Corollary~\ref{thm:trace-resistance-is-free-resistance} is a consequence of the inequality and how it tends to an equality as $k \to \iy$.
\end{remark}


\subsection{Projections in Hilbert space and the conditioning of the random walk}
\label{sec:probabilistic-interp-of-v_x-and-f_x}
  
  In Remark~\ref{rem:wired-vs-free-as-boundary-conditions}, we gave an operator-theoretic account of the difference between $R^F$ and $R^W$. The foregoing probabilistic discussions might lead one to wonder if there is a probabilistic counterpart. An alternative approach is given in \cite[App.~B]{Kig03}.
  
  On a finite network, it is well-known that 
  \begin{align}\label{eqn:vx-as-prob-on-finite}
    v_x = R(o,x) u_x,  
  \end{align}
  where $u_x(y)$ is the probability that a random walker (RW) started at $y$ reaches $x$ before $o$:
  \linenopax
  \begin{align}\label{eqn:u_x-defined-probabilistically}
    u_x(y) := \prob_y[\gt_x < \gt_o].
  \end{align}
  Here again, $\gt_x$ denotes the hitting time of $x$ as in Definition~\ref{def:hitting-time}. Note that \eqref{eqn:def:R(x,y)-energy} gives $u_x = \frac{v_x}{\energy(v_x)}$. The relationship \eqref{eqn:vx-as-prob-on-finite} is discussed in \cite{DoSn84,LevPerWil08,Lyons:ProbOnTrees}. 
  
  Theorem~\ref{thm:fx-as-prob} is a wired extension of \eqref{eqn:vx-as-prob-on-finite} to transient networks. 
  The corresponding free version appears in Conjecture~\ref{thm:v_x-as-prob-on-infinite-network}.
   
\begin{theorem}\label{thm:fx-as-prob}
  On a transient network, let $f_x$ be the representative of $\Pfin v_x$ specified by $f_x(o)=0$. Then for $x \neq o$, $f_x$ is computed probabilistically by
  \linenopax
  \begin{align}\label{eqn:prob-expr-for-fx}
    f_x(y) = R^W(o,x) \left(\prob_y[\gt_x < \gt_o] + \prob_y^{G}[\gt_o = \gt_x = \iy] \prob_x^{G}[\gt_o = \iy] \lim_{k \to \iy} \tfrac{\cond(x)}{\cond(\iy_k)} \prob_{\iy_k}^{G_k^W}[\gt_{\iy_k} < \gt_{\{x,o\}}]\right).
  \end{align}
  \begin{proof}
    Fix $x,y$ and an exhaustion $\{G_k\}_{k=1}^\iy$, and suppose without loss of generality that $o,x,y \in G_1$. 
    Since $v_x=f_x$ on any finite network, the identity \eqref{eqn:vx-as-prob-on-finite} gives $f_x^{(k)} = R_{G_k^W}(o,x) \check u_x^{(k)}$, where $f_x^{(k)}$ is the unique solution to $\Lap v = \gd_x - \gd_o$ on the finite (wired) subnetwork $G_k^W$, and 
    \linenopax
    \begin{align*}
      \check u_x^{(k)}(y) := \prob_y^{G_k^W}[\gt_x < \gt_o].
    \end{align*}
    where the superscript indicates the network in which the random walk travels.
    As in the previous case, we just need to check the limit of $\check u_x^{(k)}$, for which, we have
    \begin{align}\label{eqn:fx-as-prob-derivation-1}
       \check u_x^{(k)}(y) 
       &= \prob_y^{G_k^W}[\gt_x < \gt_o \text{ \& } \gt_x < \gt_{\iy_k}] + \prob_y^{G_k^W}[\gt_x < \gt_o \text{ \& } \gt_x > \gt_{\iy_k}]
    \end{align}
    %
    The first probability in \eqref{eqn:fx-as-prob-derivation-1} is
    \begin{align*}
      \prob_y^{G_k^W}[\gt_x < \gt_o \text{ \& } \gt_x < \gt_{\iy_k}]
      &= \prob_y^{G}[\gt_x < \gt_o \text{ \& } \gt_x < \gt_{G_k^\complm}] \\
      \limas{k}&\hstr[2.2] \prob_y^{G}[\gt_x < \gt_o \text{ \& } \gt_x < \iy] 
      = \prob_y^{G}[\gt_x < \gt_o],
    \end{align*}
    where the last equality follows because $\gt_x < \gt_o$ implies $\gt_x < \iy$.
    
    The latter probability in \eqref{eqn:fx-as-prob-derivation-1} measures the set of paths which travel from $y$ to $\iy_k$ without hitting $x$ or $o$, and then on to $x$ without passing through $o$, and hence can be rewritten
    \begin{align*}
      \prob_y^{G_k^W}&[\gt_{\iy_k} < \gt_x < \gt_o]
      = \prob_y^{G_k^W}[\gt_{\iy_k} < \gt_{\{o,x\}}] \prob_{\iy_k}^{G_k^W}[\gt_x < \gt_o] \\
      &= \prob_y^{G_k^W}[\gt_{\iy_k} < \gt_{\{o,x\}}] \left(\prob_{\iy_k}^{G_k^W}[\gt_{\iy_k} < \gt_{\{x,o\}}]\prob_{\iy_k}^{G_k^W}[\gt_x < \gt_{\{\iy_k,o\}}] + \prob_{\iy_k}^{G_k^W}[\gt_x < \gt_{\{\iy_k,o\}}]\right),
    \end{align*}
    since a walk starting at $\iy_k$ may or may not return to $\iy_k$ before reaching $x$. 
    
    First, consider only those walks which do not loop back through $\iy_k$ (i.e., multiply out the above expression and take the second term) to observe
    \begin{align}
      \prob_y^{G_k^W}[\gt_{\iy_k} < \gt_{\{o,x\}}] &\prob_{\iy_k}^{G_k^W}[\gt_x < \gt_{\{\iy_k,o\}}]
      = \prob_y^{G_k^W}[\gt_{\iy_k} < \gt_{\{o,x\}}]\prob_x^{G_k^W}[\gt_{\iy_k} < \gt_o] \tfrac{\cond(x)}{\cond(\iy_k)} \label{eqn:fx-as-prob-derivation-3} \\
      &= \left(1 - \prob_y^{G_k^W}[\gt_{\{o,x\}} < \gt_{\iy_k}]\right) \left(1 - \prob_x^{G_k^W}[\gt_o < \gt_{\iy_k}]\right) \tfrac{\cond(x)}{\cond(\iy_k)} \notag \\
      &= \left(1 - \prob_y^{G}[\gt_{\{o,x\}} < \gt_{G_k^\complm}]\right) \left(1 - \prob_x^{G}[\gt_o < \gt_{G_k^\complm}]\right) \tfrac{\cond(x)}{\cond(\iy_k)} \notag \\
      \limas{k} 
      &\hstr[2.2] \left(1 - \prob_y^{G}[\gt_{\{o,x\}} < \iy]\right) \left(1 - \prob_x^{G}[\gt_o < \iy]\right) \lim_{k \to \iy} \tfrac{\cond(x)}{\cond(\iy_k)} \notag \\
      &= \prob_y^{G}[\gt_o = \gt_x = \iy] \prob_x^{G}[\gt_o = \iy] \lim_{k \to \iy} \tfrac{\cond(x)}{\cond(\iy_k)}.
      \label{eqn:fx-as-prob-derivation-2}
    \end{align}
    Note that \eqref{eqn:fx-as-prob-derivation-3} comes by reversibility of the walk, and the way probability is computed for paths from $\iy_k$ to $x$ which avoid $o$ and $\iy_k$. Since the network is transient, $\sum_{k=1}^\iy{\cond(\iy_k)^{-1}}$ is summable by Nash-William's criterion and so $\lim_{k \to \iy} \tfrac{\cond(x)}{\cond(\iy_k)} = 0$ causes \eqref{eqn:fx-as-prob-derivation-2} to vanish. 
    
    Now for walks which do loop back through $\iy_k$, the same arguments as above yield
    \begin{align*}
      \prob_y^{G_k^W}[\gt_{\iy_k} < \gt_{\{o,x\}}] &\prob_{\iy_k}^{G_k^W}[\gt_{\iy_k} < \gt_{\{x,o\}}]\prob_{\iy_k}^{G_k^W}[\gt_x < \gt_{\{\iy_k,o\}}] \\
      &\limas{k} \hstr[1] \prob_y^{G}[\gt_o = \gt_x = \iy] \prob_x^{G}[\gt_o = \iy] \lim_{k \to \iy} \tfrac{\cond(x)}{\cond(\iy_k)} \prob_{\iy_k}^{G_k^W}[\gt_{\iy_k} < \gt_{\{x,o\}}],
    \end{align*}
    and the conclusion follows.
  \end{proof}
\end{theorem}


The following conjecture expresses a free extension of \eqref{eqn:vx-as-prob-on-finite} to infinite networks. We offer an erroneous ``proof'' in the hopes that it may inspire the reader to find a correct proof. The error is discussed in Remark~\ref{rem:error-in-prob-expr-for-vx}, just below. In the statement of Conjecture~\ref{thm:v_x-as-prob-on-infinite-network}, we use the notation
\begin{align}\label{eqn:bounded-traj}
  |\cpath| < \iy
\end{align}
to denote the event that the walk is bounded, i.e., that the trajectory is contained in a finite subnetwork of \Graph.

\begin{conj}\label{thm:v_x-as-prob-on-infinite-network}
  On an infinite resistance network, let $v_x$ be the representative of an element of the energy kernel specified by $v_x(o)=0$. Then for $x \neq o$, $v_x$ is computed probabilistically by
  \linenopax
  \begin{align}\label{eqn:prob-expr-for-vx}
    v_x(y) = R^F(o,x) \prob_y[\gt_x < \gt_o \,|\, |\cpath| < \iy],
  \end{align}
  that is, the walk is conditioned to lie entirely in some finite subnetwork as in \eqref{eqn:bounded-traj}. 
  \begin{proof}[``Proof.'']
    Fix $x,y$ and suppose without loss of generality that $o,x,y \in G_1$. One can  write \eqref{eqn:vx-as-prob-on-finite} on $G_k$ as $v_x^{(k)} = R_{G_k^F}(o,x) u_x^{(k)}$. In other words, $v_x^{(k)}$ is the unique solution to $\Lap v = \gd_x - \gd_o$ on the finite subnetwork $G_k^F$. 
    Since $R^F(x,y) = \lim_{k \to \iy} R_{G_k^F}(x,y)$ by \eqref{eqn:def:free-resistance}, it only remains to check the limit of $u_x^{(k)}$. Using a superscript to indicate the network in which the random walk travels, we have 
    \linenopax
    \begin{align}\label{eqn:prob-expr-for-vx:error}
      \lim_{k \to \iy} u_x^{(k)}(y)
      &= \lim_{k \to \iy} \prob_y^{G_k^F} [\gt_x < \gt_o]
       = \lim_{k \to \iy} \prob_y^{G} [\gt_x < \gt_o \,|\, \cpath \ci G_k^F]. 
    \end{align}
    Here again, the notation $[\cpath \ci H]$ denotes the event that the random walk never leaves the subnetwork $H$, i.e., $\gt_{H^\complm} = \iy$. The events $[\cpath \ci G_k^F]$ are nested and increasing, so the limit is the union, and \eqref{eqn:prob-expr-for-vx} follows. 
    Note that $G_k^F$ is recurrent, so $\cpath \ci G_k^F$ implies $\gt_x < \iy$.
  \end{proof}
\end{conj}

\begin{remark}\label{rem:error-in-prob-expr-for-vx}
  As indicated, the argument outlined above is incomplete due to the second equality of \eqref{eqn:prob-expr-for-vx:error}. While the set of paths from $y$ to $x$ in $G_k^F$ is the same as the set of paths from $y$ to $x$ in $G$ which lie in $G_k$, the probability of a given path may differ when computed in network or the other. This happens precisely when \cpath passes through a boundary point: the transition probability away from a point in $\bd G_k$ is strictly larger in $G_k^F$ than it is in $G_k$.
\end{remark}

\subsection{The shorted operator}
\label{sec:shorted-operator}

It is worth noting that the operator $D$ defined in \eqref{eqn:action-of-Lap-quadrants} is always invertible as in the discussion following \eqref{eqn:Schur-complement-as-sum}.
However, the Schur complement construction is valid more generally. As is pointed out in \cite{Butler}, the \emph{shorted operator} generalizes the Schur complement construction to positive operators on a (typically infinite-dimensional) Hilbert space \sH; see \cite{Anderson,AndersonTrapp,Krein}. In general, let $T=T^\ad$ be a positive operator so $\la \gf, T\gf\ra \geq 0$ for all $\gf \in \sH$, and let $S$ be a closed subspace of \sH. Partition $T$ analogously to \eqref{eqn:action-of-Lap-quadrants}, so that $A:S \to S$, $B:S \to S^\complm$, $B^T:S^\complm \to S$, and $D:S^\complm \to S^\complm$. 

\begin{theorem}[{\cite{AndersonTrapp}}]
  With respect to the usual ordering of self-adjoint operators, there exists a unique operator $\sS{h}(T)$ such that
  \begin{align*}
    \sS{h}(T) = \sup_L\left\{L \geq 0 \suth \left[\begin{array}{cc} L & 0 \\ 0 & 0\end{array}\right] \leq T\right\},
  \end{align*}
  and it is given by
  \begin{align*}
    \sS{h}(T) = \lim_{\ge \to 0^+} \left(A - B^T (D + \ge)^{-1} B\right).
  \end{align*}
  In particular, the shorted operator coincides with the Schur complement, whenever the latter exists.
\end{theorem}

There is another characterization of the shorted operator due to \cite{Butler}.

\begin{theorem}[{\cite{Butler}}]
  Suppose $\{\gy_n\} \ci \sH$ is a sequence satisfying $\la\gy_n, D\gy_n\ra \leq M$ for some $M \in \bR$, and $\lim_{n \to \iy} T \left[\begin{smallmatrix} \gf \\ \gy_n \end{smallmatrix}\right] = \left[\begin{smallmatrix} \gq \\ 0 \end{smallmatrix}\right]$. Then $\sS{h}(T) \gf = \lim_{n \to \iy} \left(A \gf + B^T\gy_n\right)$.
\end{theorem}


\section{Comparison of resistance metric to other metrics}
\label{sec:Comparison-to-other-metrics}

\subsection{Comparison to geodesic metric}
\label{sec:Comparison-to-shortest-path-metric}

On a Riemannian manifold $(\gW,g)$, the geodesic distance is
  \linenopax
  \begin{equation*}
    \distmin(x,y)
    := \inf_{\gg} \left\{\int_0^1 g(\gg'(t),\gg'(t))^{1/2} \,dt \suth \gg(0)=x, \gg(1)=y, \gg \in C^1\right\}.
  \end{equation*}
 
\begin{defn}\label{def:geodesic-metric}
  On $(\Graph, \cond)$, the \emph{geodesic distance} from $x$ to $y$ is
  \linenopax
  \begin{equation}\label{eqn:def:geodesic-metric}
    \distmin(x,y)
    := \inf \{r(\cpath) \suth \cpath \in \Paths(x,y)\},
  \end{equation}
  where $r(\cpath) := \sum_{(x,y) \in \cpath} \cond_{xy}^{-1}$. (For resistors in series, the total resistance is the sum.)
\end{defn}
  \glossary{name={$\distmin(x,y)$},description={geodesic distance, shortest-path metric},sort=d,format=textbf}
  \glossary{name={$r(\cpath)$},description={total resistance of a path},sort=r,format=textbf}

\begin{remark}\label{rem:shortest-path-metric}
  Definition~\ref{def:geodesic-metric} should not be confused with the combinatorial distance (or ``shortest path metric'') found in the literature on general graph theory. Without weights on the edges one usually defines the shortest path metric simply as the minimal number of edges in a path from $x$ to $y$. (This corresponds to taking $\cond \equiv 1$.) Such shortest paths always exist. According to Definition~\ref{def:geodesic-metric}, shortest paths may not exist (cf.~Example~\ref{exm:infinite-ladder}). Of course, even when they do exist, they are typically not unique.

It should be observed that effective resistance is not a geodesic metric, in the usual sense of metric geometry;  it does not correspond to a length structure in the sense of \cite[\S2]{Burago}. 
\end{remark}

We include the following folkloric result for completeness.

\begin{lemma}\label{thm:shortest-path-bounds-resistance-distance}
  The effective resistance is bounded above by the geodesic distance. More precisely, $R^F(x,y) \leq \distmin(x,y)$ with equality if and only if \Graph is a tree.
  \begin{proof}
    If there is a second path, then some positive amount of current will pass along it (i.e., there is a positive probability of getting to $y$ via this route). To make this precise, let $v = v_x - v_y$ and let $\cpath = (x=x_0,x_1,\dots,x_n=y)$ be any path from $x$ to $y$:
    \linenopax
    \begin{align*}
      R^F(x,y)^2
      = |v(x)-v(y)|^2
      \leq r(\cpath) \energy(v),
    \end{align*}
    by the exact same computation as in the proof of \cite[Lemma~2.5]{DGG}, but with $u=v$. The desired inequality then follows by dividing both sides by $\energy(v) = R^F(x,y)$. 
    
    The other claim follows by observing that trees are characterized by the property of having exactly one path \cpath between any $x$ and $y$ in \verts. By \eqref{eqn:def:R^F(x,y)-diss}, $R^F(x,y)$ can be found by computing the dissipation of the unit current which runs entirely along \cpath from $x$ to $y$. This means that $\curr(x_{i-1},x_i)=1$ on \cpath, and $\curr = 0$ elsewhere, so
    \linenopax
    \begin{align*}
      R^F(x,y)
      &= \diss(\curr)
      = \sum_{i=1}^n \frac1{\cond_{x_{i-1} x_i}} \curr(x_{i-1},x_i)^2 
      = \sum_{i=1}^n \frac1{\cond_{x_{i-1} x_i}}  
      = r(\cpath).
      \qedhere
    \end{align*}
  \end{proof}
\end{lemma}

\begin{remark}\label{rem:shortcut-through-infinity}
  It is clear from the end of the proof of Lemma~\ref{thm:shortest-path-bounds-resistance-distance} that on a tree, $v_x-v_y$ is locally constant on the complement of the unique path from $x$ to $y$. However, this may not hold for $f_x-f_y$, where $f_x = \Pfin v_x$; see Example~\ref{exm:geometric-integers}. This is an example of how the wired resistance can ``cheat'' by considering currents which take a shortcut through infinity; compare \eqref{eqn:def:R^F(x,y)-diss} to \eqref{eqn:def:R^W(x,y)-diss}.
\end{remark}

\subsection{Comparison to Connes' metric}
\label{sec:Comparison-to-Connes-metric}
  The formulation of $R(x,y)$ given in \eqref{eqn:def:R(x,y)-Lap} may evoke Connes' maxim that a metric can be thought of as the inverse of a Dirac operator; cf.~\cite{Con94}. This does not appear to have a literal incarnation in the current context, but we do have the inequality of Lemma~\ref{thm:Connes'-metric} in the case when $\cond = \one$. In this formulation, $v \in \HE$ is considered as a multiplication operator defined on $u$ by
  \linenopax
  \begin{equation}\label{eqn:v-as-operator}
    (vu)(x) := v(x) u(x), \q \forall x \in \verts,
  \end{equation}
  and both $v$ and \Lap are considered as operators on $\ell^2(\verts \cap \dom \energy$. We use the commutator notation $[v,\Lap] := v\Lap-\Lap v$, and $\| [v,\Lap] \|$ is understood as the usual operator norm on $\ell^2$.

\begin{lemma}\label{thm:Connes'-metric}
  If $\cond = \one$, then for all $x,y \in \verts$ one has
  \linenopax
  \begin{align}\label{eqn:Connes'-metric}
    R(x,y) \leq \sup \{|v(x)-v(y)|^2 \suth \| [v,\Lap] \| \leq \sqrt2, v \in \dom \energy\}.
  \end{align}
  \begin{proof}
    We will compare \eqref{eqn:Connes'-metric} to \eqref{eqn:def:R(x,y)-sup}. Writing $M_v$ for multiplication by $v$, it is straightforward to compute from the definitions
    \linenopax
    \begin{align*}
      (M_v \Lap - \Lap M_v)u(x)
      &= \sum_{y \nbr x} (v(y)-v(x))u(y),
    \end{align*}
    so that the Schwarz inequality gives
    \linenopax
    \begin{align*}
      \|[M_v,\Lap]u\|_2^2
      &= \sum_{x \in \verts} \left|\sum_{y \nbr x} (v(y)-v(x))u(y)\right|^2 \\
      &\leq \sum_{x \in \verts} \left(\sum_{y \nbr x} |v(y)-v(x)|^2\right) \left(\sum_{y \nbr x} |u(y)|^2\right).
    \end{align*}
    By extending the sum of $|u(x)|^2$ to all $x \in \verts$, this (admittedly crude) estimate gives \linebreak \mbox{$\|[v,\Lap]u\|_2^2 \leq 2 \|u\|_2^2 \energy(v)$,} and hence $\|[v,\Lap]\|^2 \leq 2 \energy(v)$.
  \end{proof}
\end{lemma}

\subsection{Effective resistance between measures}
\label{sec:Generalized resistance metrics}
\label{sec:Effective-resistance-between-measures}
We describe a notion of effective resistance between probability measures, of which $R(x,y)$ ($R^F$ or $R^W$) is a special case. This concept is closely related to the notion of total variation of measures, and hence is related to mixing times of Markov chains; cf. \cite[\S4.1]{LevPerWil08}. When the Markov chain is taken to be random walk on an ERN, the state space is just the vertices of \Graph.

\begin{defn}\label{def:resistance-distance-of-measures}
  Let \gm and \gn be two probability measures on \verts. Then the total variation distance between them is
  \linenopax
  \begin{align}\label{eqn:def:resistance-distance-of-measures}
    \dist_{\TV}(\gm,\gn)
    := 2 \sup _{A \ci \verts} |\gm(A)-\gn(A)|.
  \end{align}
\end{defn}

\begin{prop}[{\cite[Prop.~4.5]{LevPerWil08}}]
  \label{thm:resistance-distance-of-measures-as-norm}
  Let \gm and \gn be two probability measures on the state space \gW of a (discrete) Markov chain. The total variation distance between them is
  \linenopax
  \begin{align}\label{eqn:resistance-distance-of-measures-as-norm}
    \dist_{\TV}(\gm,\gn)
    = \sup\left\{ \left|\sum_{x \in \gW} u(x) \gm(x) - \sum_{x \in \gW} u(x) \gn(x)\right| \suth \|u\|_\iy \leq 1 \right\}.
  \end{align}
  Here, $\|u\|_\iy := \sup_{x \in \verts}|u(x)|$.
\end{prop}
  \glossary{name={$\dist_{\TV}(\gm,\gn)$},description={total variation metric},sort=d,format=textbf}

If we think of \gm as a linear functional acting on the space of bounded functions, then it is clear that \eqref{eqn:resistance-distance-of-measures-as-norm} expresses $\dist_{\TV}(\gm,\gn)$ as the operator norm $\|\gm-\gn\|$. That is, it expresses the pairing between $\gm \in \ell^1$ and $u \in \ell^\iy$. We can therefore extend $R^F$ directly; see Remark~\ref{rem:R(x,y) = |L{xy}|}.

\begin{defn}
  The \emph{free effective resistance} between two probability measures \gm and \gn is
  \linenopax
  \begin{align}
    \dist_{R^F}(\gm,\gn) 
    := \sup\left\{ \left|\sum_{x \in \verts} u(x) \gm(x) - \sum_{x \in \verts} u(x) \gn(x)\right|^2 \suth \|u\|_\energy \leq 1 \right\},
  \end{align} 
  in accordance with \eqref{eqn:def:R^F(x,y)-sup}, and the \emph{wired effective resistance} between them is
  \linenopax
  \begin{align}
    \phantom{MMM}\dist_{R^W}(\gm,\gn) 
    := \sup\left\{ \left|\sum_{x \in \verts} u(x) \gm(x) - \sum_{x \in \verts} u(x) \gn(x)\right|^2 \suth \|u\|_\energy \leq 1, u \in \Fin \right\},
  \end{align} 
  in accordance with \eqref{eqn:def:R^W(x,y)-sup}.
\end{defn}
  \glossary{name={$\dist_{R^F}(\gm,\gn)$},description={free resistance between two probability measures},sort=d,format=textbf}

It is clear from this definition (and Remark~\ref{rem:R(x,y) = |L{xy}|}) that $R^F(x,y) = \dist_{R^F}(\gd_x,\gd_y)$ and $R^W(x,y) = \dist_{R^W}(\gd_x,\gd_y)$. This extension of effective resistance to measures was motivated by a question of Marc Rieffel in \cite{Rieffel99}.

\section{von Neumann construction of the energy space \HE}
\label{sec:Construction-of-HE}
\label{sec:vonNeumann's-embedding-thm}

In Theorem~\ref{thm:R^F-embed-ERN-in-Hilbert} we show that an \ERN equipped with resistance metric may be embedded in a Hilbert space in such a way that $R^{1/2}$ is the norm difference of the corresponding vectors in the Hilbert space. It turns out that (up to unitary isomorphism) the Hilbert space is \HE when the embedding is applied with $R=R^F$ and \Fin when applied with $R=R^W$; see Remark~\ref{rem:uniqueness-of-HE}. As a consequence, we obtain an alternative and independent construction of the Hilbert space \HE of finite-energy functions. This provides further justification for \HE as the natural Hilbert space for studying the metric space $(\Graph,R^F)$ and \Fin as the natural Hilbert space for studying the metric space $(\Graph,R^W)$. We use the notation $(\Graph, R)$ when the distinction between $R^F$ and $R^W$ is not important. 

\begin{defn}\label{def:negative-semidefinite}
  A function $d:X \times X \to \bR$ is \emph{negative semidefinite} iff for any $f:X \to \bR$ satisfying $\sum_{x \in X} f(x) = 0$, one has
  \linenopax
  \begin{equation}\label{eqn:negative-semidefinite}
    \sum_{x,y \in F} f(x) d^2(x,y) f(y) \leq 0,
  \end{equation}
  where $F$ is any finite subset of $X$.
\end{defn}

\begin{theorem}[von Neumann]\label{thm:vonNeumann's-embedding-thm}
  Suppose $(X,d)$ is a metric space. There exists a Hilbert space $\sH$ and an embedding $w:(X,d) \to \sH$ sending $x \mapsto w_x$ and satisfying
  \linenopax
  \begin{equation}\label{eqn:vonNeumann's-embedding-thm}
    d(x,y) = \|w_x - w_y\|_\sH
  \end{equation}
  if and only if $d^2$ is negative semidefinite. 
\end{theorem}

von Neumann's theorem also has a form of uniqueness which may be thought of as a universal property.

\begin{theorem}[von Neumann]\label{thm:uniqueness-of-vNeu-embedding}
  If there is another Hilbert space \sK and an embedding $k:\sH \to \sK$, with $\|k_x-k_y\|_\sK = d(x,y)$ and $\{k_x\}_{x \in X}$ dense in \sK, then there exists a unique unitary isomorphism $U:\sH \to \sK$.
\end{theorem}

The following theorem is inspired by the work of von Neumann and Schoenberg \cite{vN32a,Ber96,Schoe38a,Schoe38b,Ber84} on ``screw functions'', but is a completely new result. One aspect of this result that contrasts sharply with the classical theory is that the embedding is applied to the metric $R^{1/2}$ instead of $R$, for each of $R=R^F$ and $R=R^W$.

\begin{theorem}\label{thm:R^F-embed-ERN-in-Hilbert}
  $(\Graph,R^F)$ may be isometrically embedded in a Hilbert space.
  \begin{proof}
    According to Theorem~\ref{thm:vonNeumann's-embedding-thm}, we need only to check that $R^F$ is negative semidefinite (see Definition~\ref{def:negative-semidefinite}). Let $f:\verts \to \bR$ \saty $\sum_{x \in \verts} f(x) = 0$. We must show
that $\sum_{x,y \in F} \cj{f(x)} R^F(x,y) f(y) \leq 0$, for any finite subset $F \ci \verts$. From \eqref{eqn:def:R^F(x,y)-energy}, we have
    \linenopax
    \begin{align*}
      \sum_{x,y \in F} \cj{f(x)} R^F(x,y) f(y)
      &= \sum_{x,y \in F} \cj{f(x)} \energy(v_x - v_y) f(y) \\
      &= \sum_{x,y \in F} \cj{f(x)} \energy(v_x) f(y) - 2 \sum_{x,y \in F} \cj{f(x)} \la v_x, v_y \ra_\energy f(y) + \sum_{x,y \in F} \cj{f(x)} \energy(v_x) f(y) \\
      &= -2\left\la \sum_{x \in F} f(x) v_x, \sum_{y \in F} f(y) v_y \right\ra_\energy \\
      &= -2\left\|\sum_{x \in F} f(x) v_x \right\|_\energy^2
      \leq 0.
    \end{align*}
    For the second equality, note that the first two sums vanish by the assumption on $f$.
  \end{proof}
\end{theorem}

\begin{cor}\label{thm:R^W-embed-ERN-in-Hilbert}
  $(\Graph,R^W)$ may be isometrically embedded in a Hilbert space.
  \begin{proof}
    Because the energy-minimizer in \eqref{eqn:def:R^W(x,y)-energy} is $f_x = \Pfin v_x$, we can simply repeat the proof of Theorem~\ref{thm:R^F-embed-ERN-in-Hilbert} with $f_x$ in place of $v_x$.
  \end{proof}
\end{cor}

\begin{cor}\label{rem:uniqueness-of-HE}
  Up to unitary isomorphism, the von Neumann embedding of $(\Graph,R^F)$ into a Hilbert space is $V:\Graph \to \HE$ given by $x \mapsto v_x$, and the von Neumann embedding of $(\Graph,R^W)$ into a Hilbert space is $\sF:\Graph \to \Fin$ given by $x \mapsto f_x$. 
  \begin{proof}
    Since $R^F(x,y) =\|v_x-v_y\|_\energy^2$ by \eqref{eqn:def:R^F(x,y)-energy}, Theorem~\ref{thm:uniqueness-of-vNeu-embedding} shows that the embedded image of $(\Graph,R^F)$ is unitarily equivalent to the \energy-closure of $\spn\{v_x\}$, which is \HE. Similarly, $R^W(x,y) =\|f_x-f_y\|_\energy^2$, where $f_x := \Pfin v_x$, by \eqref{eqn:def:R^W(x,y)-energy}, whence the embedded image of $(\Graph,R^W)$ is unitarily equivalent to the \energy-closure of $\spn\{f_x\}$, which is \Fin.
  \end{proof}
\end{cor}

von Neumann's theorem is constructive, and provides a method for obtaining the embedding, which we briefly describe, continuing in the notation of Theorem~\ref{thm:vonNeumann's-embedding-thm}.

\begin{FlatList}
  \item Schwarz inequality. If $d$ is a negative semidefinite function on $X \times X$, then define a positive semidefinite bilinear form on functions $f,g:X \to \bC$ by
      \linenopax
      \begin{align} \label{eqn:vNeu-inner-prod}
        Q(f,g) = \la f,g\ra_Q := - \sum_{x,y} f(x) d^2(x,y) g(y).
      \end{align}
      This gives a quadratic form $Q(f) := Q(f,f)$ for which 
      $Q(f,g)^2 \leq Q(f) Q(g)$ holds.
  \item The kernel of $Q$. Denote the collection of finitely supported functions on $X$ by $\Fin(X)$ and define $\Fin_0(X) := \{f \in \Fin(X) \suth {\textstyle \sum}_x f(x) = 0\}$. Before completing $\Fin_0(X)$ with respect to $Q$, 
  define the subspace $\ker Q = \{f \in \Fin_0(X) \suth Q(f)=0\}$.
  \item Pass to quotient. Define $\tilde Q$ to be the induced quadratic form on the quotient space $\Fin_0(X)/\ker Q$. Now $\tilde Q$ is \emph{strictly} positive definite on the quotient space and $\|\gf\|_{\sH_{vN}} := -\tilde Q(\gf)$ will be a bona fide norm.
  \item Complete. Completing the quotient space with respect to $\|\cdot\|_{\sH_{vN}}$ gives a Hilbert space
      \linenopax
      \begin{align}\label{eqn:vNeu's-quotient-completion}
        \sH_{vN} &:= \left(\frac{\Fin_0(X)}{\ker Q}\right)^\sim,
        \q\text{with}\q
        \la \gf, \gy \ra_{\sH_{vN}} = -\tilde Q(\gf,\gy),
      \end{align}
      into which $(X,d)$ may be embedded.
\end{FlatList}

\begin{remark}\label{rem:norm-vs-quasinorm}
  In the construction outlined above, one can choose any vertex $o \in \verts$ to act as the ``origin'' and it becomes the origin of the new Hilbert space $\sH_{vN}$. As a quadratic form defined on the space of all functions $v:\verts \to \bC$, the energy is indefinite and hence allows one to define only a quasinorm. One way to deal with the fact that \energy does not ``see constant functions'' is to adjust the energy so as to obtain a true norm:
  \begin{equation}\label{eqn:energy-o}
    \la u,v \ra_o := u(o)v(o) +  \la u,v \ra_\energy.
  \end{equation}
This is done in \cite{Yamasaki79,Kayano88,MuYaYo,Kig01,Lyons:ProbOnTrees}, for example, and a comparison with the current approach is discussed in \cite[\S4.1]{DGG}.
  We have instead elected to work ``modulo constants'' because 
  the kernel of \energy is the set of constant functions, and inspection of von Neumann's embedding theorem (cf.~\eqref{eqn:vNeu's-quotient-completion}) shows that it is precisely these functions which are ``modded out'' in von Neumann's construction.
\end{remark}

\begin{cor}\label{thm:H_energy-as-functions-of-finite-distance}
  $v \in \HE$ if and only if $\sum_{x,y \in \verts} v(x) R(x,y) v(y) < \iy$.
\end{cor}

Corollary~\ref{thm:H_energy-as-functions-of-finite-distance} follows from the construction outlined above 
  \version{}{\marginpar{Where does this come from? \\ Palle says: uniqueness in vN-Sch}
  \[\sum_{x,y \in \verts} v(x) R(x,y) v(y) = \left\|\sum_{x \in \verts} v(x)(\gd_x-\gd_o)\right\|_\energy^2,\]
  }
and is comparable to $G(f,f)$ in \cite[\S2]{Kayano82}.

\section{Examples}
\label{sec:examples}

In this section, we introduce the most basic family of examples that illustrate our technical results and exhibit the properties (and support the types of functions) that we have discussed above. 


  \begin{figure}
    \centering
    \scalebox{1.0}{\includegraphics{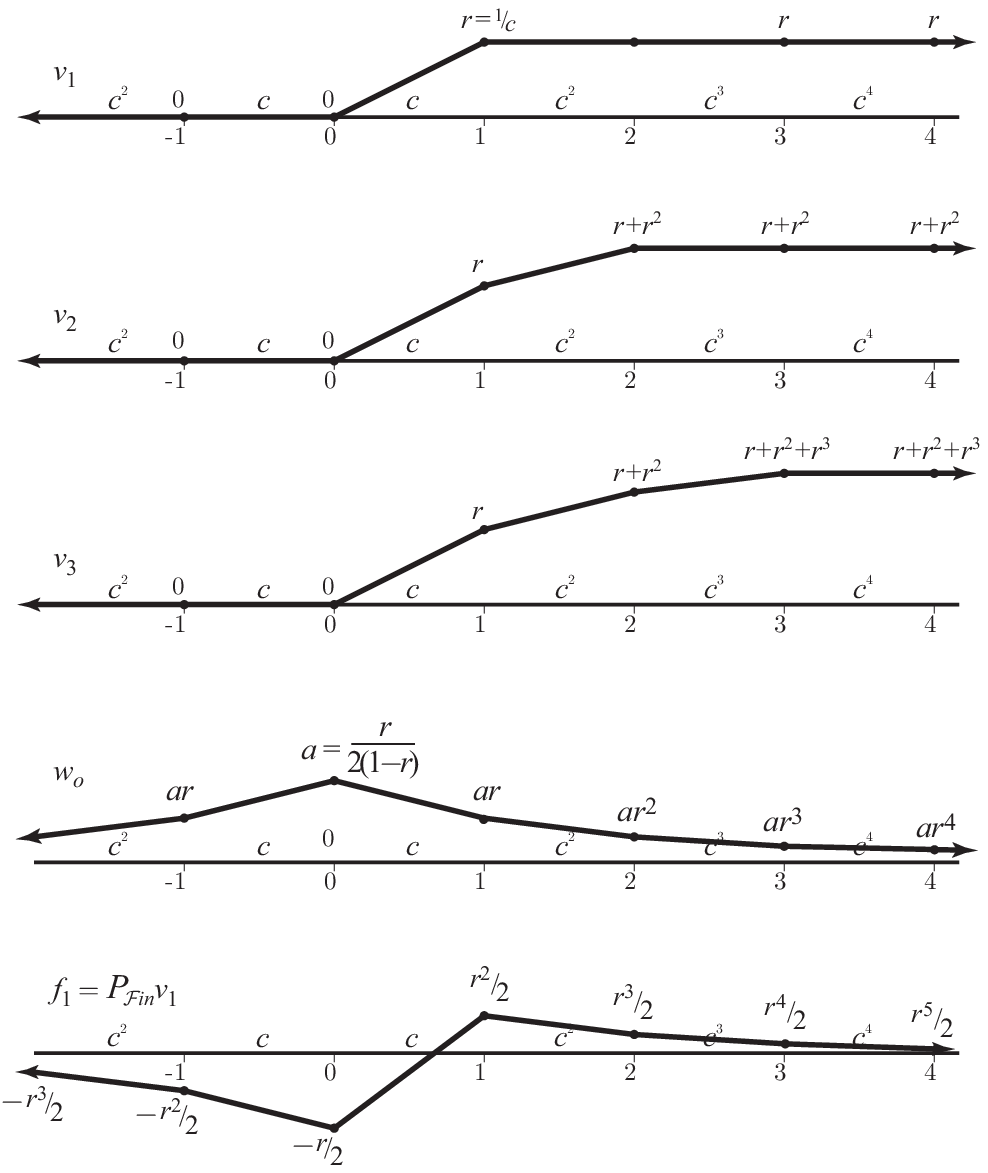}}
    \caption{\captionsize The functions $v_1$, $v_2$, and $v_3$ on $(\bZ,c^n)$. Also, the monopole $w_o$ and the projection $f_1 = \Pfin v_1$. See Example~\ref{exm:geometric-integers}}
    \label{fig:vx-on-Zc}
  \end{figure}

\begin{exm}[Integer networks]\label{def:integers}
  Let $(\bZ,\cond)$ denote the network with integers for vertices, and with conductances defined by \cond.
  We fix $o=0$.
\end{exm}

These networks are more interesting when \cond tends to \iy. For example, $\sum \cond_{xy}^{-1} < \iy$ implies $\Harm \neq 0$ for $(\bZ,\cond)$, as is shown in \cite[Thm.~6.4]{DGG}. It is helpful to keep the following more concrete model in mind, especially if one hopes for tractable computations.

\begin{exm}[Geometric integer model]\label{exm:geometric-integers}
  For a fixed constant $c>1$, let $(\bZ,c^n)$ denote the network with integers for vertices, and with geometrically increasing conductances defined by $\cond_{n-1,n} = c^{\max\{|n|,|n-1|\}}$ so that the network under consideration is
  \linenopax
  \begin{align*}
    \xymatrix{
      \dots \ar@{-}[r]^{c^3}
      & -2 \ar@{-}[r]^{c^2} 
      & -1 \ar@{-}[r]^{c} 
      & 0 \ar@{-}[r]^{c} 
      & 1 \ar@{-}[r]^{c^2} 
      & 2 \ar@{-}[r]^{c^3} 
      & 3 \ar@{-}[r]^{c^4} 
      & \dots
    }
  \end{align*} 
  Again, we fix $o=0$. On this network the energy kernel is given by
  \linenopax
  \begin{align*}
    v_n(k) = 
    \begin{cases}
      0, &k \leq 0, \\
      \frac{1-r^{k+1}}{1-r}, &1 \leq k \leq n, \\
      \frac{1-r^{n+1}}{1-r}, &k \geq n,
    \end{cases}
    n > 0,
  \end{align*}
  and similarly for $n < 0$.
  Also, the function $w_o(n) = ar^{|n|}$, $a:= \frac{r}{2(1-r)}$, defines a \emph{monopole} (that is, $w_o \in \HE$ and $\Lap w_o = \gd_o$), and $h(n) = \operatorname{sgn}(n) (1-w_o(n))$ defines an element of \Harm. See Figure~\ref{fig:vx-on-Zc}.
  
  In Figure~\ref{fig:vx-on-Zc}, one can also see that $f_1 = \Pfin v_1$ induces a current flow of 1 amp from 1 to 0, with  $\frac{1+r}2$ amps flowing down the 1-edge path from 1 to 0, and the remaining current of $\frac{1-r}2$ amps flowing down the ``pseudo-path'' from 1 to $+\iy$ and then from $-\iy$ to 0.
\end{exm}
  
\begin{exm}[Binary tree]
  \label{exm:Binary-tree}
  One may have $\Harm \neq 0$ for networks with constant conductances $\cond = \one$, provided they branch sufficiently rapidly. For example, consider the binary tree $(\sT,\one)$. An exhaustion of this network is depicted in Figure~\ref{fig:tree-exhaustions} and Figure~\ref{fig:tree-repkernels} depicts an element of the energy kernel $v_x$ and its projections to \Fin and \Harm. 
  Let $o$ be the root, and let $V_n$ be the set of vertices of $(\sT,\one)$ that are $n$ steps from $o$. If one were to solder these vertices $V_n$ together, for each $n$, the result would be isomorphic to $(\bZ,c^n)$ for $c=2$.\footnote{Nonetheless, there remain important differences between $(\sT,\one)$ and $(\bZ,2^n)$. For example, \Lap is self-adjoint on the former but not on the latter, due to the presence of nontrivial deficiency spaces on $(\bZ,2^n)$ for $c>1$; see \cite[Ex.~13.41]{OTERN}.}
  \begin{figure}
  \centering
  \scalebox{1.0}{\includegraphics{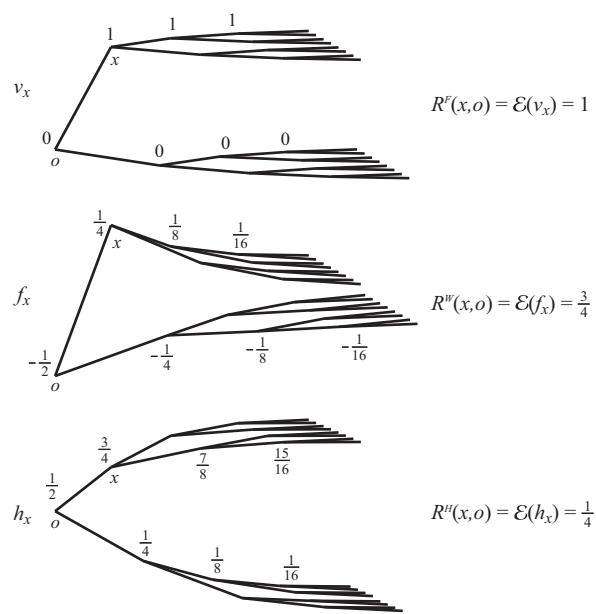}}
  \caption{\captionsize The reproducing kernel on the tree with $\cond = \one$. For a vertex $x$ which is adjacent to the origin $o$, this figure illustrates the elements $v_x$, $f_x = \Pfin v_x$, and $h_x = \Phar v_x$; see Example~\ref{exm:Binary-tree}. }
  \label{fig:tree-repkernels}
\end{figure}
\end{exm}
  
\begin{exm}[Infinite ladder]
  \label{exm:infinite-ladder}
  The following infinite one-sided ladder model furnishes a situation where no shortest path exists, as mentioned in Remark~\ref{rem:shortest-path-metric}. Resistances are as labelled.
  \begin{equation}\label{eqn:exm:one-sided-ladder-model-no-shortest-path}
    \xymatrix{
    *+[l]{x \; \bullet} \ar@{-}[r]_{\frac14} \ar@{-}[d]_1
      & \bullet \ar@{-}[r]_{\frac1{16}} \ar@{-}[d]_{\frac1{4}} & \bullet \ar@{-}[r]_{\frac1{64}} \ar@{-}[d]_{\frac1{16}} & \bullet \ar@{-}[r] \ar@{-}[d]_{\frac1{64}}
      &\dots \ar@{-}[r]_{4^{-n}} & \bullet \ar@{-}[d]_{4^{-n}} \ar@{-}[r]_{4^{-(n+1)}} & \dots \\
    *+[l]{y \; \bullet} & \bullet \ar@{-}[l]_{\frac14} & \bullet \ar@{-}[l]_{\frac1{16}} & y_3 \ar@{-}[l]_{\frac1{64}} &\dots  \ar@{-}[l] & \bullet \ar@{-}[l]_{4^{-n}} \ar@{-}[r]^{4^{-(n+1)}} & \dots
    }
  \end{equation}
  Then $\dist_\cpath(\ga,\gw) = \frac23$, but $r(\cpath) > \frac23$ for every path \cpath from $x$ to $y$.
\end{exm}

\bibliographystyle{alpha}
\bibliography{ERM}

\end{document}